\documentclass[letter,fleqn,12pt]{article} 
\usepackage{graphicx}
\usepackage{subfigure}
\usepackage[T1]{fontenc}
\usepackage{mathrsfs, calligra, pbsi, frcursive, la, aurical}

\usepackage{bbold, epsfig}
\usepackage{latexsym, amsmath, amssymb, euscript, eufrak}

\usepackage{fontenc, indentfirst}
\usepackage{fancyhdr}
\usepackage{multicol}
\usepackage{slashbox}
\usepackage{stmaryrd}
\usepackage{aurical}
\usepackage[usenames,dvipsnames]{xcolor}
\usepackage[english]{babel}%french
\usepackage{ucs}
\usepackage[latin1]{inputenc}

\usepackage[numbers]{natbib}

\reversemarginpar

\usepackage{amsmath}
\usepackage{amssymb}
\usepackage{theorem}
\usepackage{ifthen}

\newtheorem{theorem}{Theorem}[section] 
 
\newtheorem{proposition}[theorem]{Proposition} 
 
\newtheorem{remark}[theorem]{Remark}

\newcommand{\sgn}[1]{{\rm sgn}{#1}}

\newcommand{\ol}[1]{\overline{#1}}
\newcommand{\mR}{\mathbb R}
\newcommand{\ds}[1]{\displaystyle{#1}}

\newcommand{\colre}[1]{\textcolor{black}{#1}}
\newcommand{\colbl}[1]{\textcolor{black}{#1}}

\begin{document}

\begin{center}
\textbf{\Large 
Modeling, shape analysis and computation of the equilibrium pore shape near a PEM-PEM intersection} 
\end{center}

\begin{center}
 Peter {\sc Berg}\footnote{NTNU, Norway}, 
 Sven-Joachim {\sc Kimmerle}\footnote%[3]
{Universit\"at der Bundeswehr M\"unchen, Germany}, 
 Arian {\sc Novruzi}\footnote{University of Ottawa, Canada}$^{,}$%$^{2,}$
\footnote[4]{Corresponding author}
\end{center}

\begin{figure}[!h]
\begin{center}
\includegraphics[width=80mm,height=40mm]{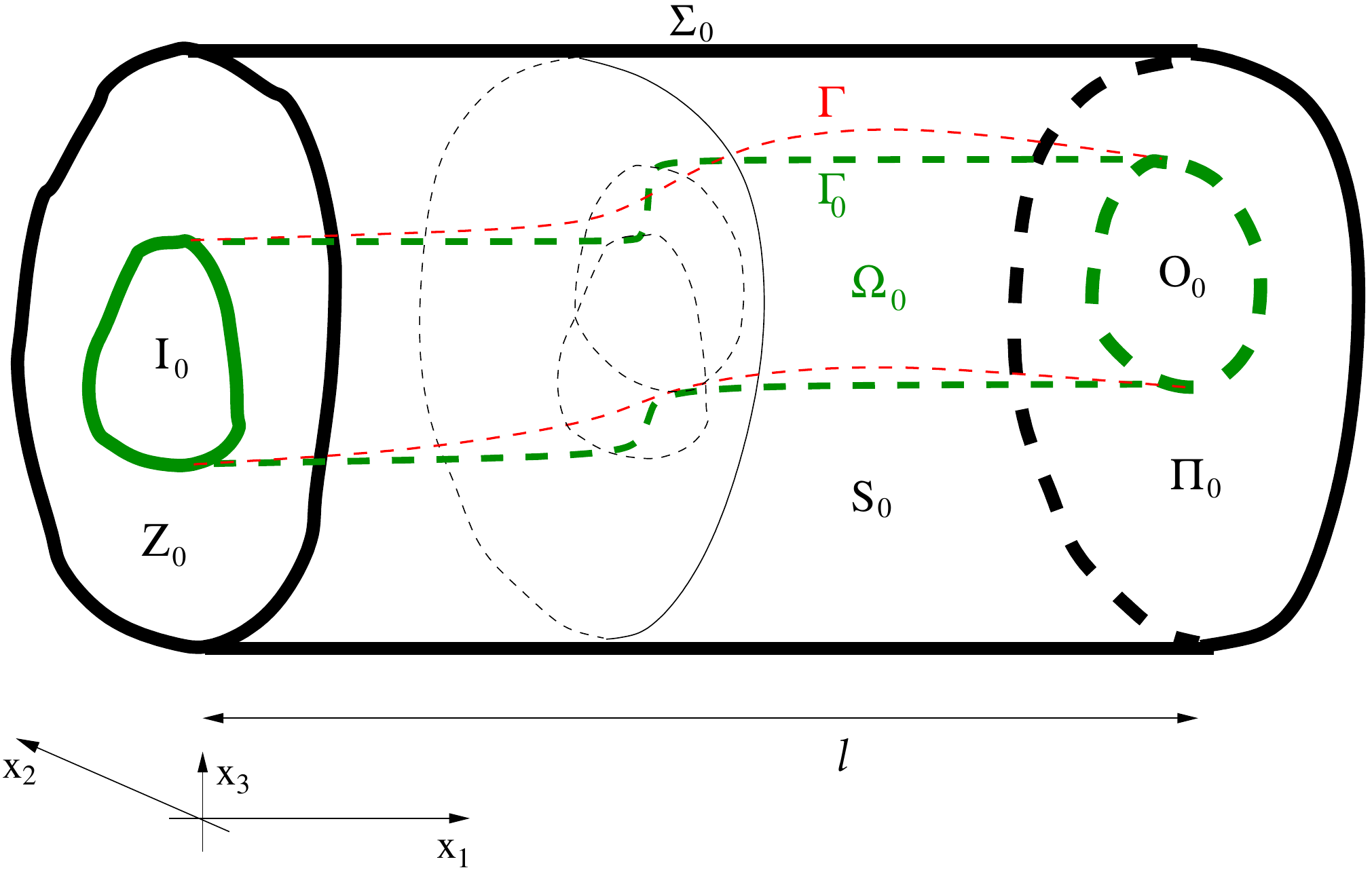}
\end{center}
\caption{\label{f:D+O}The domain $\Omega_0$ \colre{(the union of two nano-channels)}, $S_0$ \colre{(the union of two elastomers)} and the boundaries}
\end{figure}

{\small\noindent
{\bf Abstract:}
\colre{In this paper we study} the equilibrium shape of an interface that represents the \colre{lateral} boundary of a pore channel embedded in an elastomer. The model consists of a system of PDEs, comprising a linear elasticity equation for displacements within the elastomer and a nonlinear Poisson equation for the electric potential within the channel (filled with protons and water). To determine the equilibrium interface, a variational approach is employed. We analyze: i) the existence and uniqueness of the electrical potential, ii) the shape derivatives of state variables and iii) the shape differentiability of the corresponding energy and the corresponding Euler-Lagrange equation. The latter leads to a modified Young-Laplace equation on the interface. This modified equation is compared with the classical Young-Laplace equation by computing several equilibrium shapes, \colre{using} a fixed point algorithm.


\section{Introduction}
In this contribution, we study the equilibrium shape \colre{of } an interface which represents
the \colre{lateral} \colbl{boundary} of a pore channel embedded in an elastomer (solid elastic body). This problem
originates from the modeling of the electrical resistance between two adjacent polymer electrolyte membranes, a material used in hydrogen fuel cells.

In mathematical terms, we consider a system of partial differential equations (PDEs) which is a simplified version of the elasticity, Stokes and Nernst-Planck equations in the absence of velocity and any external electric field.

%paragraph by SJ about some background of PEM and fuel cells
The motivation for considering the problem under investigation is to understand the \colre{interaction} between charged fluid flow and the morphology of the fluid domain, i.e. the interface between the fluid and the elastomer, in a polymer electrolyte membrane (PEM).
PEMs are an essential component of the so-called proton exchange membrane fuel cells. 
This type of fuel cell, running at low temperature, converts hydrogen and oxygen into electric energy, and is expected to power automobiles in the not-too-distant future. 
Within the PEM fuel cell, hydrogen enters the device at the (negative) anode and is oxidized \colre{at the anode catalyst layer}, producing protons and electrons. The protons migrate across the PEM, a charge-selective
medium, to the cathode, \colre{and} the electrons flow through an outer circuit \colre{to reach the cathode}. 
Meanwhile, the oxygen enters the fuel cell at the (positive) cathode and enters into a reaction with the hydrogen 
protons and the electrons that arrive at the cathode catalyst layer. 
As a result, useful electric current is produced, with water as a byproduct.

Polymer electrolyte membranes are made from ionomers, which consist of long hydrophobic 
backbones with shorter hydrophilic side chains \cite{kreuer}, \cite{mauritz}. 
The latter end in acid groups, such as sulfonic 
acid groups in Nafion, which enter into an ion (proton) exchange equilibrium when the PEM 
adsorbs water. Minimization of the system's free energy leads to phase separation where 
proton-conducting water pores form and include the acid groups, surrounded by hydrophobic 
domains. Two key aspects of PEM research are i) the proton and water flow inside the pores, and 
ii) the pore formation itself, i.e. the membrane morphology, related to water uptake. 

It is widely believed that PEMs consist of many small nanochannels of cylinder-like shape \cite{SchmidtRohrChen}. Since the surface dynamics of the PEM will determine PEM 
functionality inside a fuel cell to a large degree, it remains an interesting question how two PEMs, 
and their nanoscopic pores, will connect across %the
 their mutual surface when pressed together, see e.g.~the discussion about stack design without bipolar plates in \cite{RecentTrendsCh3}.
It is known that an interface resistance arises, the cause of which, however, is unexplained.
 
Hence, we focus on two such nanochannels filled with water and protons, that (partially) connect, as it may occur at the interface of two PEMs. In Fig.~\ref{f:D+O}, the nanochannel and the elastomer corresponding to one PEM, are represented by one half of $\Omega_0\cup S_0$ \colre{(the subscript/superscript $0$ is used for the domain and variables related to the reference configuration)}.

The ohmic resistance between $I_0$ and $O_0$ is of particular importance.
To address this question, one must first study the (equilibrium) shape of the elastomer/fluid flow interface $\Gamma$, and how its shape depends on the type of connection between two nanochannels. Note that the equilibrium interface $\Gamma$ minimizes an energy, whose Euler-Lagrange equation leads to a modified Young-Laplace equation.

For sake of comparing our modified Young-Laplace equation with the literature results, we present also 
the classical Young-Laplace equation and a fixed point algorithm associated with it, used for solving the equilibrium interface.
We present several numerical examples, which demonstrate the differences and similarities between our modified Young-Laplace equation and the classical Young-Laplace equation, and give some conclusions.

For a model of this problem in the case of a radially symmetric channel $\Omega_0$ and neglecting the elastomer $S_0$, see \cite{LBKN-11}, where the effects of parameters on proton transport in  nanopores are analyzed numerically. 
A more general model for a PEM pore taking into account charged fluid flow, external electric field and fluid structure interaction, is presented and examined numerically (by a fixed point approach) in \cite{KBN-12}.

\section{Mathematical formulation}
Let $D_0\subset \mathbb R^N$, $N=2,3$,  be a nonempty simply connected, open, bounded and fixed domain, $\Omega_0\varsubsetneq D_0$ a nonempty simply connected, open set and $S_0=D_0\backslash\overline{\Omega}_0$. 
\colre{The boundaries of these domains are denoted as follows: 
$\partial\Omega_0=I_0\cup \Gamma_0\cup O_0$,
$\partial S_0=Z_0\cup\Sigma_0\cup\Gamma_0\cup \Pi_0$},
$\partial D_0:=I_0\cup Z_0\cup \Sigma_0\cup \Pi_0 \cup O_0$}.
We assume \colre{that $I_0\neq\emptyset$, $O_0\neq\emptyset$, 
$I_0\cup Z_0\subset\{x_1=0\}$, $O_0\cup \Pi_0\subset\{x_1=\ell\}$},  see Fig. \ref{f:D+O}.
Furthermore, let $\nu_{\colre{0}}$ be the normal vector to $\partial D_0$ or $\partial\Omega_0$, exterior to 
$D_{\colre{0}}$ or $\Omega_0$ \colre{(note that the subscript/superscript $0$ is related to the reference configuration and the associated variables)}. 
%and 
%$F_\tau^0=F^0 - (F^0\cdot\nu_0)\nu_0$ denote the tangential component of any vector $F^0$ on $\partial D_0$ or $\Gamma_0$.

Here, the domain $\Omega_0$ represents the \colre{initial} space occupied by an electrically charged fluid and $S_0$ represents the space occupied by the elastomer.

In $S_0$, the unknown variable is the displacement $U^0=(U_1^0,U_2^0,U_3^0)\in H^1(S_0,\mR^N)$
satisfying
\begin{eqnarray}
-\nabla \cdot \sigma^0(U^0)
&=&
0\quad in\ S_0. 	\label{e:U0}
\end{eqnarray}
Here $\sigma^0(U^0)$ is the first Piola-Kirchhoff stress tensor given by
\begin{eqnarray}
\sigma^0(U^0)&=& (\sigma_{ij}^0(U^0))
=
\left(-p_S^0\delta_{ij}+\epsilon^0_{ij}(U^0)\right)		\nonumber\\
&:=&
-
p_S^0\delta_{ij}
+
\left(k_S - \frac{2}{3} G_S\right) (\nabla U^0)_{ij} \delta_{ij} + 
 G_S ((\nabla U^0)_{ij} + (\nabla U^0)_{ji}),	 		 	\label{e:sigma}
\end{eqnarray}
with 
\colre{$p_S^0$ a given solid reference pressure,
$k_S>0$ and $G_S>0$, the bulk modulus of the elastic material and its shear modulus, respectively 
(see Table \ref{t:data}), satisfying
\begin{equation}\label{e:kSGS}
k_S-\frac{2}{3}G_S\geq0. 
\end{equation}
}
Equation (\ref{e:U0}) is equipped with the following boundary conditions
\begin{eqnarray}
U^0
&=&
0\quad on\ Z_0\cup\Sigma_0\cup\Pi_0, \quad \mbox{\it(no-slip boundary condition)}			\label{bc:U0-Z0}\\
%
%\sigma^0_\tau(U^0)=0,\quad
%U^0\cdot{\nu_0}
%&=&
%0
%\quad
%on\ \Sigma_0,\quad (slip\ on\; \Sigma_0),			\label{bc:U0-S0}\\
%\sigma^0_\tau(U^0)=0,\quad
%U^0\cdot{\nu_0} 
%&=&
%\lambda^0\cdot\nu_0=\lambda^0_1
%\quad on\; \Pi_0,								\label{bc:U0-P0}\\
U^0
&=&
\lambda\quad on\ \Gamma_0.			    		\label{bc:U0-G0}
\end{eqnarray}
%where $\sigma^0_\tau(U^0)=\sigma^0(U^0)\cdot\nu_0-(\nu_0\cdot\sigma^0(U^0)\cdot \nu_0)\nu_0$ is the tangential stress %and $\lambda^0$ is a known displacement. 
The displacement $\lambda$ on $\Gamma_0$ is unknown and is defined by an equilibrium force balance equation, see section \ref{s:nm}. Note that $\lambda$ must satisfy some compatibility conditions on
$\partial \Gamma_0$ so that the electric potential is smooth enough, see sections 
\ref{sss:min-pb},  \ref{sss:existence}, \ref{sss:shd-u}.

Set $\mathcal{F}=\mathcal{F}(U^0):=I+U^0$, where $I$ is the identity transformation, $I(X)=X$, $X\in S_0$.
Under the transformation $\mathcal{F}$, the domain $S_0$ is transformed to the domain 
$S_\lambda=\mathcal F(S_0):=\{\mathcal{F}(X),\ X\in S_0\}$ \colre{(the actual (deformed) elastomer configuration)}. 
\colre{Note that the compatibility conditions on $\lambda$, see section \ref{sss:min-pb}, ensure $S_\lambda\subset D_0$
and $\partial S_\lambda=Z_0\cup\Sigma_0\cup\Gamma_\lambda\cup O_0$, where $\Gamma_\lambda=\mathcal F(\Gamma_0)$.
Then we set $\Omega_\lambda=D_0\backslash S_\lambda$ \colre{(the actual (defoemed) channel configuration)}, so we have $\partial\Omega_\lambda=I_0\cup\Gamma_\lambda\cup O_0$, see Fig. \ref{f:D+O}}.
Furthermore, let $\nu_\lambda$  be the normal vector on $\Gamma_\lambda$, exterior to $\Omega_\lambda$. 

\colre{We emphasize that the subscript/superscript $\lambda$ is related to the actual (deformed) configuration.
Wherever permissible, we will drop the subscript/superscript $\lambda$}.

In $\Omega$ the unknown variables are $\varphi$ (the electric potential, primary unknown), $c$ (the proton concentration) and the pressure $p$ (both $c$ and $p$, the secondary unknowns depending on $\varphi$). The variable $\varphi$ satisfies
\begin{eqnarray}
-\nabla\cdot(\varepsilon\nabla \varphi) 
&=&
Fc \quad in\  \Omega,			\label{e:phi}\\
%\int_\Omega\varphi &=&0,\\
-\varepsilon\partial_\nu\varphi 
&=&
\sigma_c\chi_\Gamma 
\quad on\ \partial\Omega,		\label{bc:phi}
\end{eqnarray}
with  $\chi_\Gamma$ the characteristic function of $\Gamma$,
$\sigma_c>0$, $- \sigma_c$ representing the surface charge density of negative sulfonic acid groups, and
$\varepsilon = \varepsilon_0 \varepsilon_r$, the electric permittivity in vacuum. 
\colre{Here $\varepsilon_0= 8.8542 \times 10^{-12}$ F/m is a universal constant}
and $\varepsilon$ is a parameter, assumed constant  \colre{(see Table \ref{t:data} for the values of parameters)}.

%%\footnote{normally, the surface charge of acid groups is negative, say $\sigma<0$, 
%%and the boundary condition for $\varphi$ is 
%%$-\epsilon\partial_\nu\varphi =-\sigma\chi_\Gamma$; we denote $\sigma$ instead of $-\sigma$ and then  $\sigma>0$}.

Note that from (\ref{e:phi}), (\ref{bc:phi}) it follows that $\varphi$ satisfies the following necessary condition
\begin{eqnarray}
F\int_{\Omega} c 
&=&
\int_{\Gamma} \sigma_c
\quad \mbox{\it(global electro-neutrality (GEN))} .	\label{e:gen}
\end{eqnarray}
%%%referred as ``global electro-neutrality condition''.
%
%For the uniqueness of $\varphi$ we impose also the condition
%\begin{equation}\label{e:Iphi=0}
% \int_I \varphi = 0.
%\end{equation}
%So, for $\varphi$ smooth there are points in $I$ such that $\varphi=0$.

The variables $c$ and $p$ satisfy the following equations \colbl{(\cite{LBKN-11}, equations~(11), (12))}
\begin{eqnarray}
\nabla c+ \frac{F}{RT} c\nabla \varphi 
&=&
0\quad in\  \Omega,	\label{e:c}\\
\nabla p + F c \nabla \varphi 
&=&
0\quad in\ \Omega. 	\label{e:p}
\end{eqnarray}
These equations are derived from the stationary Nernst-Planck and Stokes equations, respectively, when we neglect the pressure drop along the channel (so the velocity is zero) and the external electric field.
\colre{From (\ref{e:c}), (\ref{e:p})}  follows 
\begin{equation}\label{e:c,p=}
 \ds{c=c_0e^{-\frac{F}{RT}\varphi}},\quad
 p=RT(c-c_0)+p_0.
\end{equation}
Here, $p_0$ and $c_0$ are given and independent of $\lambda$. 
The constant $p_0$ represents the pressure value associated with the location where the concentration is $c_0$, and hence $\varphi=0$.

We emphasize that the displacement $\lambda$ is unknown. It is defined such that the interface $\Gamma=\mathcal F(\Gamma)$ is at equilibrium, i.e. a force balance on $\Gamma$ is achieved. 
The objective of this paper is to identify the boundary equations  \colre{which determine} this interface $\Gamma$, and to compute $\Gamma$ numerically for different physical parameters.

To this end, we will consider a variational method. 
\colre{The associated functional represents the energy of the system, namely the mechanical, elastic and surface tension energies}.
Associated with this method, we will discuss the existence of the solution to (\ref{e:phi}), (\ref{bc:phi}), which is a challenging problem and deserves attention in its own right.
We will derive the formula for the shape derivative of the energy affiliated with our problem, which 
leads to a free boundary equation referred to as a modified Young-Laplace equation. We will conclude with some numerical results of the interface~$\Gamma$.

To place our variational approach within the context of known results, we will restate the classical Young-Laplace equation.
We will compute numerically the interface $\Gamma$ based on  this equation and a fixed point method 
\colre{(like in \cite{quarteroni-1})}, and compare these results with those obtained with our modified Young-Laplace equation.

\section{Free boundary equation}\label{s:nm}
In this section, we describe our variational approach to the equilibrium interface $\Gamma$.
This approach is based on the minimization of a certain (Lagrangian) energy. The Euler-Lagrange 
equation associated with the minimization problem related to this energy provides a new force balance equation (a modified Young-Laplace equation), which leads to a gradient descent minimization algorithm.

In order to relate our method to known results, we will recall first the commonly used classical Young-Laplace  equation which represents a direct balance of elastic, hydrodynamic and surface tension forces.
This Y-L equation leads to a fixed point numerical algorithm.

\subsection{Classical Young-Laplace equation and fixed point method}
It is accepted in the literature, see for example \cite{Wilmanski-08}, that the interface $\Gamma$ is at equilibrium if 
the elastic forces, the \colre{hydrodynamical (pressure)} and the surface tension forces are balanced, i.e.
\begin{equation}\label{e:fb}
\sigma\cdot\nu 
+
p\nu
=
\gamma \mathcal H \nu
\quad
on\ \Gamma.
\end{equation}
Here, $\mathcal H$ is the mean curvature \colre{(here, it means the sum of two principal curvatures)} of $\Gamma$ seen from $\Omega$, $\sigma=\sigma(U)$ is 
the Cauchy stress tensor in $S$ (the actual configuration of the solid) and $\gamma$ is the surface tension 
on $\Gamma$.
Note that (\ref{e:fb}) represents the classical Young-Laplace equation.

Here, it is assumed that on $\Gamma$ only solid stresses, pressure and surface tension forces act.
At this point we emphasize that when considering the equilibrium interface from a variational viewpoint (see Section \ref{ss:vm}), it turns out that  \colre{additional forces \colbl{act} on $\Gamma$}.

The Cauchy stress tensor $\sigma$ (in $S=(I+U^0)(S_0)$) is related to the 1st Piola-Kirchhoff stress tensor 
$\sigma^0$ (in $S_0$) by (see \cite{beatty-1, ciarlet-1, ciarlet-2})
\begin{eqnarray}
\sigma\circ \mathcal{F}
&=&
(\sigma^0\cdot[\nabla \mathcal{F}]^t)|det[\nabla \mathcal{F}]^{-1}|,		\label{e:sigma=sigma0}\\
(\sigma\cdot\nu) d\Gamma 
&=&
(\sigma^0\cdot{\nu_0})d\Gamma_0,					\label{e:sigma.n=sigma0.n0}\\
d\Gamma
&=&
|det[\nabla \mathcal{F}][\nabla \mathcal{F}]\cdot{\nu_0}|d\Gamma_0.			\label{e:dG=dG0}
\end{eqnarray}

By using (\ref{e:sigma=sigma0}), (\ref{e:sigma.n=sigma0.n0}) and (\ref{e:dG=dG0}), we can write equation (\ref{e:fb}) equivalently on $\Gamma_0$ as follows
\begin{equation}\label{e:fb0}
\sigma^0(U^0)\cdot{\nu_0}
=
\left((\gamma\mathcal H - p)\nu\right)\circ \mathcal{F}\ |det[\nabla \mathcal{F}][\nabla \mathcal{F}]\cdot{\nu_0}|
\quad
on\ \Gamma_0.
\end{equation}

To solve (\ref{e:fb}), or equivalently (\ref{e:fb0}), we use a fixed point argument as in \cite{quarteroni-1}. 
Namely, let us consider the so-called Dirichlet-to-Neumann operator defined by
\begin{equation}\label{e:D-N}
\begin{array}{rccc}
A:&H^{1/2}(\Gamma_0,\mR^N)&\to& H^{-1/2}(\Gamma_0,\mR^N), \\
  &\lambda& \to&A(\lambda)=\sigma^0\cdot{\nu_0}, 
\end{array}
\end{equation}
where $\sigma^0=\sigma^0(U^0)$ and $U^0$ is the solution of (\ref{e:U0}), (\ref{bc:U0-Z0})-(\ref{bc:U0-G0%P0
}). 
Note that $A$ is well defined (see Theorem \ref{th:U}). Then, (\ref{e:fb0}) is equivalent to
\[
A(\lambda)
= 
((\gamma\mathcal H-p)\nu)\circ \mathcal{F}\ |det[\nabla \mathcal{F}][\nabla \mathcal{F}]\cdot{\nu_0}|
\quad\ on\ \Gamma_0. 
\]
Let $A^{-1}:H^{-1/2}(\Gamma_0,\mR^N)\mapsto H^{1/2}(\Gamma_0,\mR^N)$ be the inverse of $A$ (as in Theorem \ref{th:U}, one can easily prove that $A^{-1}$ is well-defined). Applying $A^{-1}$ to both sides of the last equation gives
\begin{equation}\label{e:fp}
 \lambda 
 = A^{-1}(((\gamma\mathcal H-p)\nu)\circ \mathcal{F}|det[\nabla \mathcal{F}][\nabla \mathcal{F}]\cdot\nu_0|) 
 =:B(\lambda)\quad on\ \Gamma_0.
\end{equation}
So, $\lambda\in H^{1/2}(\Gamma_0,\mR^N)$ is a fixed point of the operator $B$. This implies the following fixed point algorithm 
for solving (\ref{e:fb})
\begin{equation}\label{e:fpm}
\left\{
{\rm 
\begin{tabular}{rl}
Step 0.& $n=0$. Initialize $\lambda$. \\
Step 1.& a) Compute $U^0$ in $S_0$.\\
       & b) Compute $\varphi$, $c$, $p$ in $\Omega=D_0\backslash\ol{S}$, $S=(I+U^0)(S_0)$.\\	
       & c) Compute $B(\lambda)$ and set $\lambda=B(\lambda)$.\\
       & d) Set $n=n+1$.\\
Step 2.& Repeat Step 1 until $\|\lambda-B(\lambda)\|_{H^{1/2}(\Gamma_0;\mathbb R^N)} \leq err$,\\
\end{tabular}
}
\right.
\end{equation}
where $err$ is a fixed given error.

We will not present any analysis related to this approach. We will only use the algorithm (\ref{e:fpm}) for computing 
numerically the equilibrium interface, so we can compare it with the numerical results obtained with our variational
method \colre{(see \colbl{section} \ref{s:nr})}.

\subsection{Variational method}\label{ss:vm}
For $\lambda\in\Lambda$, with $\Lambda$ a set of admissible deformations to be specified later,
see section \ref{sss:min-pb},  let $U^0$ be defined \colre{in $S_0$} by (\ref{e:U0}), (\ref{bc:U0-Z0}), (\ref{bc:U0-G0}),
$\varphi$ be defined  in \colre{$\Omega_\lambda$} by (\ref{e:phi})-(\ref{bc:phi}),
and  $c$ and $p$ be given by (\ref{e:c,p=}). 

We consider the energy functional $E(\lambda)$ given by
\begin{eqnarray}
 E(\lambda)
 &=&
 E_{mech} + E_{el} + E_{st},			\label{e:E}
\end{eqnarray}
with 
\begin{eqnarray}
\hspace*{-5mm}
E_{mech}
&:=&
\int_{S_0}
\left(\frac{1}{2}\epsilon^0(U^0):[\nabla U^0]-p_S^0\nabla\cdot U^0 \right)
-
\int_\Omega p
=:E_{mech,s} + E_{mech,l},						\label{e:Emech}\\
\hspace*{-5mm}
E_{el}
&:=&
-
\left(
\int_{\Omega} \frac{\varepsilon}{2}|\nabla\varphi|^2+\int_\Gamma\sigma_c\varphi
\right),									\label{e:Eel}\\
\hspace*{-5mm}
E_{st}
&:=&
\int_\Gamma \gamma,							\label{e:Est}\\
\hspace*{-5mm}
\epsilon^0(U^0):[\nabla U^0]
&:=&
\epsilon^0_{ij}(U^0)\partial_i U^0_j,				\nonumber
\end{eqnarray}
and $p=RT(c-c_0)+p_0$, $\ds{c=c_0e^{-\frac{F}{RT}\varphi}}$, as given by (\ref{e:c,p=}).
We have written $\epsilon^0(U^0)$ to emphasize the dependence of $\epsilon^0$ on $U^0$   
(we will use similar notations for $\sigma^0(U^0)$).

\begin{remark}\label{r:E}
Note that if we define 
\begin{eqnarray*}
e(V^0,\psi,\lambda)
&=&
\int_{S_0}
\frac{1}{2}\epsilon^0(V^0):[\nabla V^0]-p_S^0\nabla\cdot V^0 
-
\int_{\Omega_\lambda} RTc_0(e^{-\frac{F}{RT}\psi}-1)+p_0	\\
&-&
\left(
\int_{\Omega_\lambda} \frac{\varepsilon}{2}|\nabla\psi|^2+\int_{\Gamma_\lambda}\sigma_c\psi
\right)
+
\int_{\Gamma_\lambda} \gamma,
\end{eqnarray*}
then $E(\lambda)=e(U^0,\varphi,\lambda)$.
If $\lambda$ is fixed, and $e$ is considered as a function of $(V^0,\psi)\in H^1(S_0,\mR^N)\times H^1(\Omega)$,
then the first variation of $e$ with respect to $(V^0,\psi)$ gives (\ref{e:U0}), 
%(\ref{bc:U0-Z0}), (\ref{bc:U0-G0}),  
(\ref{e:phi}) and (\ref{bc:phi}).
Indeed, for example for  (\ref{e:phi}) and (\ref{bc:phi}), if $\psi\in C^\infty(\overline{\Omega})$, by differentiating $e$ at $(V^0,\varphi)$ with respect to $\psi$ gives
\begin{eqnarray*}
\partial_{\psi} e(V^0,\varphi,\lambda)(\psi)
&=&
\int_{\Omega_\lambda}
-\varepsilon(\nabla\varphi\cdot\nabla\psi)
+
Fc_0e^{-\frac{F}{RT}\varphi}\psi 
-
\int_{\Gamma_\lambda}\sigma_c\psi		\\
&=&
\int_{\Omega_\lambda}
(\nabla\cdot(\varepsilon\nabla\varphi)+Fc)\psi 
-
\int_{\partial{\Omega_\lambda}}(\varepsilon\partial_\nu\varphi + \sigma_c\chi_{\Gamma_\lambda})\psi.
\end{eqnarray*}
Therefore, if $\partial_{\psi}e(V^0,\varphi,\lambda)(\psi)=0$ \colre{for any $\psi\in H^1(\Omega)$}, 
then from the arbitrariness of $\psi$ we get (\ref{e:phi}), (\ref{bc:phi}). 

Similarly, if $\partial_{V^0}e(U^0,\psi,\lambda)(V^0)=0$ for any $V^0\in H^1_0(S_0,\mR^N)$ and 
$U^0$ satisfies the boundary conditions (\ref{bc:U0-Z0})-(\ref{bc:U0-G0}), then $U^0$ satisfies  (\ref{e:U0}).

Note also that from (\ref{e:phi}), (\ref{bc:phi}) we have
\[
\int_\Gamma \sigma_c\varphi +  \int_\Omega \varepsilon|\nabla\varphi|^2 = \int_\Omega Fc\varphi,
\]
which, by using (\ref{e:c,p=}), gives
\begin{eqnarray}
E_{el}
&=&
\int_{\Omega} 
\frac{\varepsilon}{2}|\nabla\varphi|^2
-
Fc\varphi
=
\int_\Omega
\frac{\varepsilon}{2}|\nabla\varphi|^2
+
RTc \ln\frac{c}{c_0}.					\label{e:Eel+} 
\end{eqnarray}
The last term of $E_{el}$ represents the entropic contribution of the protons. Therefore, 
$E$ should be regarded as the free energy of the system rather than its internal energy. 
Using the above expression \colre{for $E_{el}$}, we obtain this equivalent form of $E$
\begin{equation}\label{e:E+}
\hspace*{-5mm}
E(\lambda)
=
\int_{S_0}
\frac{1}{2}\epsilon^0(U^0):[\nabla U^0]-p_S^0\nabla\cdot U^0
+
\int_{\Omega} 
\left(
\frac{\varepsilon}{2}|\nabla\varphi|^2
+
RTc \ln\frac{c}{c_0}
- p
\right)
+
\int_\Gamma \gamma.		
\end{equation}
\end{remark}

\subsubsection{Minimization problem and main result}\label{sss:min-pb}
\colre{
For the initial interface $\Gamma_0$ we assume also that it meets $\partial D_0$ at an angle $\pi/2$. Namely, we assume
\begin{equation}\label{e:gamma0}
\Gamma_0=\{x\in D_0,\;\, \gamma_0(x_1,x_2,x_3)=0\},\;\;
\gamma_0\mbox{\it\;\;  is\; }  C^2,\quad 
 \nabla\gamma_0\cdot e_1=0\;\; on\; \partial \Gamma_0,
\end{equation}
where $\cdot$ is the inner product and $e_1=(1,0,0)$.
}
%%This assumption will ensure appropriate regularity for $\varphi$, see sections \ref{sss:min-pb}, \ref{sss:existence}, \ref{sss:shd-u}.
The space $\Lambda$ of admissible $\lambda$ is defined as 
\begin{eqnarray}
\Lambda
&=&\{\lambda=(\lambda_1,\ldots,\lambda_N)\in C^2(\ol{D}_0, \mR^N),\; 	\nonumber\\
&&
\;\;
\lambda=0\;\; on\;\, \partial D_0,			\nonumber\\
&&
\;\;
\partial_1\lambda_2\partial_2\gamma_0 + \partial_1\lambda_3\partial_3\gamma_0=0\;\, on\;\; \partial\Gamma_0
\}.			\label{e:L}
\end{eqnarray}

\begin{remark}
\colre{
First we note that $\Lambda\subset C^2(\ol{D}_0,\mR^N)$ and equipped with the norm of $C^2(\ol{D}_0,\mR^N)$ is a Banach space.
Note also that  if $\|\lambda\|_{C^2(\ol{D}_0,\mR^N)}$ is small, then $I+\lambda$ is invertible.
Therefore, the condition $\lambda=0$ on $\partial D_0$ ensures that $(I+\lambda)(D^0)=D_0$.
}

\colre{
The last condition on (\ref{e:L}), which is  called ``compatibility condition'', ensures that $\Gamma$ ($=\Gamma_\lambda$) meets 
$\partial D_0$ at an angle $\pi/2$. This condition follows easily because from 
$\Gamma_0=\{x\in D_0,\; \gamma_0(x)=0\}$ we obtain  
$\Gamma=\{y\in D_0,\; \gamma(y):=\gamma_0\circ(I+\lambda)^{-1}(y)=0\}$. 
Then the condition $\nabla\gamma\cdot e_1=0$  on $\partial\Gamma$ , which ensures that $\Gamma$ meets $\partial D_0$ at an angle $\pi/2$,  together with (\ref{e:gamma0}) give the last condition of (\ref{e:L}).
This condition will be used when proving the existence and the regularity of the electrical potential $\varphi$ (see theorems \ref{th:u-exist+uniq}, \ref{th:u'}).
}
\end{remark}
We look for a solution of the problem:
\begin{equation}\label{e:min(E)}
\mbox{\it find $\lambda\in \Lambda$ such that }\quad
E(\lambda)=\min\{E(\mu),\; \mu\in\Lambda\}.
\end{equation}

If $\lambda$ is a solution of (\ref{e:min(E)}), then $\Gamma_\lambda$ is the interface where the forces are in balance.
We may find the Euler-Lagrange equation associated with (\ref{e:min(E)}), see Subsection \ref{sss:shd-E}. 
\begin{proposition}\label{p:E'=}
\colre{Assume (\ref{e:gamma0}) holds and} let $\lambda\in\Lambda$ such that $I+\lambda$ is invertible. Then:\\
i)
The  functional $\mu\in\Lambda\mapsto E(\mu)$ is differentiable at $\lambda$ from $\Lambda$
to $\mathbb R$  and for $\mu\in\Lambda$ we have
\begin{eqnarray}
\hspace*{-5mm}
\frac{d}{d\mu}E(\lambda)(\mu)
=
\int_{\Gamma}
 \left(
 \Big(
 -\sigma\cdot\nu
 -\left(p+\frac{\varepsilon}{2}(|\partial_\tau\varphi|^2-(\sigma_c/\varepsilon)^2)\right)\nu
 \Big) 
+
 (\gamma - \sigma_c\varphi)\mathcal H \nu
 \right)
  \cdot\mu.
 &&					\label{e:E'(l)(m)}
\end{eqnarray}
ii) If $\lambda$ is a solution of (\ref{e:min(E)}), then 
\begin{equation}\label{e:fb-v}
\sigma\cdot\nu 
+ 
\Big(p+\frac{\varepsilon}{2}(|\partial_\tau\varphi|^2-(\sigma_c/\varepsilon)^2)\Big)\nu
=
\left(\gamma - \sigma_c\varphi\right)\mathcal H  \nu
\quad
on\ \Gamma,
\end{equation}
where $\partial_\tau$ is the tangential operator defined by $\partial_\tau f=\nabla f - (\nabla f\cdot\nu)\nu$ on $\Gamma$.
\end{proposition}

\begin{remark}\label{r:extra-term}
%Comparing (\ref{e:fb-v}) and (\ref{e:fb})  we see that they are different.
If we set
\[
p_* = p + \frac{\varepsilon}{2}(|\partial_\tau\varphi|^2-(\sigma_c/\varepsilon)^2),\qquad
\gamma_* = \gamma - \sigma_c\varphi,
\]
then  (\ref{e:fb-v}) is written as
\begin{equation}\label{e:fb-v+}
\sigma\cdot\nu 
+
p_*\nu
=
\gamma_*\mathcal H \nu \;\; on\;\, \Gamma.
\end{equation}
Comparing (\ref{e:fb}) and (\ref{e:fb-v+}), the equation (\ref{e:fb-v+}) indicates that due to the mechanical and electrical energy terms in $E(\lambda)$, the ``pressure force`` on $\Gamma$ is $p_*$ rather than  $p$, and that the surface tension %force
 is $\gamma_*$ rather than $\gamma$.

Let us point out that in the absence of the electrical charges we have $\sigma_c=c=0$ and then  (\ref{e:fb-v+}) is equivalent to the classical Young-Laplace equation (\ref{e:fb}).
\end{remark}

The equation (\ref{e:fb-v+}) implies the following  algorithm for solving $\lambda$ numerically,
or equivalently the equilibrium interface $\Gamma$, solution of (\ref{e:min(E)}).
\begin{equation}\label{e:gm}
\hspace{-6mm}
\left\{
\begin{tabular}{rp{120mm}}
Step 0.& $n=0$. Initialize $\lambda$. \\
Step 1.& a) Compute $U^0$ in $S_0$. \\
       & b) Compute $\varphi$, $c$, $p$ in $\Omega=D_0\backslash\ol{S}$, where $S=(I+U^0)(S_0)$ .\\	
       & c) Compute  
	$g= -\sigma\cdot\nu  - p_*\nu  + \gamma_* \mathcal H \nu$ 
      and set  
       $\lambda= \lambda - kg$, where $k>0$ is a small appropriate  parameter.\\
       & d) Set $n=n+1$.\\
Step 2.& Repeat Step 1 until $|E'(\lambda)|\leq err$.
\\
\end{tabular}
\right.
\end{equation}

\subsubsection{Existence and uniqueness of the state variables}\label{sss:existence}
We will prove that the energy functional $E$ is differentiable with respect to $\lambda$ and we will find the formula for its derivative.
\colre{To this end we will first prove that the state variables $U^0$ and $\varphi$ are well defined}.

\subsubsection*{The existence and uniqueness of the displacement $U^0$}
Note that any $\lambda\in\Lambda$ satisfies the boundary conditions  (\ref{bc:U0-Z0})-(\ref{bc:U0-G0}).
\begin{theorem}\label{th:U}
Let $\lambda\in\Lambda$.
Then the problem (\ref{e:U0}), (\ref{bc:U0-Z0})-(\ref{bc:U0-G0}) has a unique weak solution 
$U^0\in H^1(S_0,\mR^N)$, 
\colre{which can be written in the form $U^0=\lambda+V^0$, with $V^0\in H^1_0(S_0,\mR^N)$ the unique solution of (\ref{e:V0-weak})}.
\end{theorem}
{\bf Proof}.
\colre{The proof of the existence and uniqueness of a weak solution $U^0$ is classical, even under slightly weaker condition that (\ref{e:kSGS}, see Theorem 6.3-5, \cite{ciarlet-1} (or \cite{ciarlet-2}). For simplicity, let us present the proof}.

We look for $U^0=\lambda+V^0$. From (\ref{bc:U0-Z0})-(\ref{bc:U0-G0}) it follows that $V^0\in H^1_0(S_0,\mR^N)$.
\colre{Multiplying the $i$-th equation of (\ref{e:U0}) by $W^0_i$, where $W^0=(W^0_1,\ldots,W^0_N)\in H^1_0(S_0,\mR^N)$,
adding them, integrating by parts  and using (\ref{bc:U0-Z0})-(\ref{bc:U0-G0}) leads to
\begin{eqnarray}
b(V^0,W^0)
:=
\int_{\colre{S_0}}  
\epsilon^0(V^0):\nabla W^0  
&=&
\int_{\colre{S_0}} (p_S^0 I - \epsilon^0(\lambda)):\nabla W^0
=:
l(W^0),					\label{e:V0-weak}
\end{eqnarray}
for all $W^0\in H^1_0(S_0,\mR^N)$, where
$a:b=\sum_{i,j=1,N}a_{ij}b_{ij}$ for all $a=(a_{ij})$, $b=b_{ij}$ in $\mR^{N\times N}$.
}

\colre{
Note that $l\in H^{-1}(S_0,\mR^N)$ and as $\epsilon^0(V^0)$ is symmetric, it follows that
\[
 b(V^0,W^0)=\int_{S_0} 
 \left(k_S-\frac{2}{3}G_S\right)\sum_{i=1,N}\partial_i V_i^0\partial_i W^0_i 
 +
 \frac{G_S}{2}
 \sum_{i,j=1,N}(\partial_j V^0_i + \partial_iV^0_j)(\partial_j W^0_i + \partial_iW^0_j).
\]
Hence, $b(V^0,W^0)$ is a bilinear symmetric continuous and coercive form in $H^1_0(S_0,\mR^N)$ (the coerciveness 
follows from Korn's inequality (see Theorem 6.3-3, \cite{ciarlet-1} and (\ref{e:kSGS})).
Then the existence and uniqueness of a solution $V^0\in H^1_0(S_0,\mR^N)$ to (\ref{e:V0-weak}) 
follows trivially from Lax-Milgram lemma in $H^1_0(S_0,\mR^N)$. 
}
\hfill$\Box$

\colre{
\begin{remark}\label{r:U0-classical}
 Note that the right hand side $l(W^0)$can be written as
 \[
  l(W^0)=\int_{S_0}f(x)W^0(x)dx,\;\;\;
  f(x)=\nabla\cdot\sigma^0(\lambda).
 \]
As $\lambda\in C^2(\overline{S}_0,\mR^N)$ it follows $f\in C^0(\overline{S}_0,\mR^N)$,
hence $f\in  L^q(S_0,\mR^N)$, $q\geq 6/5$. Then from theorem 6.3-6, \cite{ciarlet-1}, it follows
$V^0\in W^{2,q}(S_0,\mR^N)$, so $U^0\in W^{2,q}(S_0,\mR^N)$ for all $q>1$. From Sobolev embeddings theorems it follows
$U^0\in C^{1,1-\delta}(\overline{S}_0,\mR^N)$,  $\delta\in(0,1)$ (by taking $q$ large) and therefore $U^0$ is a classical solution.
\end{remark}
}

\subsubsection*{The existence and uniqueness of the electrical potential $\varphi$}
Let us first note that if we set $\ds{u=-\frac{F}{RT}}\varphi$, then $u$ solves
\begin{eqnarray}
 -\Delta u + u_0 e^u &=& 0\quad in\;\, \Omega,\quad with\; u_0=c_0\frac{F^2}{RT\varepsilon},	\label{e:u}\\
 \partial_\nu u &=& g\chi_\Gamma\quad on\;\, \partial\Omega,\quad g = \sigma_c\frac{F}{RT\varepsilon}.	\label{e:bcu}
\end{eqnarray}

We look for a weak solution \colre{of (\ref{e:u}), (\ref{e:bcu}), given by}
\begin{eqnarray}
\mbox{\it find }\  u\in H^1(\Omega),\quad
\int_\Omega (\nabla u\cdot\nabla v) + u_0(e^u v) &=& \int_\Gamma g v,\quad	\forall v\in H^1(\Omega).  \label{e:uw}
\end{eqnarray}
Without restriction we may assume $u_0=1$. In fact, if \colre{$\hat{u}$} is the solution of (\ref{e:uw}) for $u_0=1$, then
$u=\colre{\hat{u}}-\ln u_0$ solves (\ref{e:uw}) for arbitrary $u_0>0$.

The equation $-\Delta u + e^u=0$ has been considered in literature in a more general context, namely in the form 
$-\Delta u + g(u)=\mu$, with $g$ increasing, $g(0)=0$, and $\mu$ a measure.

In the case of homogeneous Dirichlet boundary conditions\colre{:}
\colre{in} \cite{brezis-1} it is proven that for $\mu\in L^1(\Omega)$ the problem has a unique solution;
in \cite{brezis-2} it is proven that this problem does not have a
solution in general - more precisely, in \cite{brezis-2} the ``good'' measures $\mu$ are studied, i.e. those for which the problem has a solution, as well as the properties of these measures;
in \cite{bartolucci-1} it is proven that $-\Delta u + e^u-1=\mu$ has a unique solution in $H^1_0(\Omega)$ if $\mu\leq 4\pi H^{N-2}$, where $H^{N-2}$ denotes the $N-2$ dimensional Hausdorff measure $N\geq 3$.

In the case of homogeneous Neumann boundary conditions\colre{:}
\colre{in} \cite{hunlich-1} the equation 
$-\Delta u + e(u)=f$ is considered with an appropriate $e(u)$ (essentially $e'\geq\gamma>0$), for which the uniqueness of the solution is proven by straightforward arguments;
in \cite{hiranao-1} it is proven that $-\Delta u - \lambda u + a e^u=\epsilon f$ has at least two $H^1(\Omega)$ solutions for a range of $(\lambda,\epsilon,a)$ - the proof is made by approximating the solution in the space of eigenfunctions of $-\Delta$ with Neumann boundary conditions.

The difficulty of solving the problem (\ref{e:u}), (\ref{e:bcu}) is due to the nonlinearity $e^u$ and the Neumann 
boundary condition. In fact one can consider a variational solution to (\ref{e:u}), (\ref{e:bcu}) as in \cite{brezis-2}.
But due to the lack of Poincar\'e's inequality in $H^1(\Omega)$, which is related to the Neumann boundary condition, the minimizing sequence of the associated variational problem a priori is not bounded.
We will show that, in fact, the minimizing sequence is bounded in $H^1(\Omega)$. This will allow to extract a converging subsequence, which will provide a $W^{2,q}(\Omega)$ solution to (\ref{e:uw}).
More precisely, we will prove the following result.
\begin{theorem}\label{th:u-exist+uniq}
\colre{Assume (\ref{e:gamma0}) and} let $\lambda\in\Lambda$ with $I+\lambda$ invertible and set $\Omega=(I+\lambda)(\Omega_0)$.
Then for any $q\in(1,\infty)$, there exists a unique solution $u\in W^{2,q}(\Omega)$  of (\ref{e:uw}).
\end{theorem}
{\bf Proof}.
The proof is made in several steps.\\
i) 
Consider the functional $G:H^1(\Omega)\to\mathbb R$ defined by
$\ds{G(u)=\int_\Omega\frac{1}{2}|\nabla u|^2  + e^u - \int_\Gamma gu}$,
and look for a solution $u\in H^1(\Omega)$ of
\begin{equation}\label{e:min-G}
 G(u)=\min\left\{G(v),\; v\in H^1(\Omega),\; \int_\Omega e^u <\infty \right\}.
\end{equation}
Let $u_n$, $n\in\mathbb N$ be a minimizing sequence of $G$ in $H^1(\Omega)$. 
Then {\it $u_n$ is bounded in $H^1(\Omega)$}.
Indeed, let $u_n=u_n^++u_n^-$, where $u_n^{\pm}$, is the positive/negative part of $u_n$ (so, $u_n^+\geq0$,
$u_n^-\leq0$). It is well-known (see for example \cite{gilbarg+trudinger-1}) that $u_n^+, u_n^-\in H^1(\Omega)$.
From $G(u_n)\leq G(u_n^-)$  we get
\begin{equation}\label{e:|u+|H1<}
 \int_\Omega|\nabla u_n^+|^2 + e^{u_n^+} 
 \leq
 C\left(|\Omega| + \int_\Gamma g u_n^+\right).
\end{equation}
As $e^{u_n^+}\geq \frac{1}{2}|u_n^+|^2$, (\ref{e:|u+|H1<}) implies
\begin{equation}\label{e:|un+|H1<}
 \int_\Omega|\nabla u_n^+|^2 + |u_n^+|^2 
 \leq
 C
 \left(
 \int_\Gamma |g| |u_n^+| 
 +
 |\Omega|\right)
 \leq 
 C(\Omega,g) + \frac{1}{2}\|u_n^+\|_{H^1(\Omega)}^2,
\end{equation}
which implies that $u_n^+$ is bounded in $H^1(\Omega)$.

From $G(u_n)\leq G(u_n^+)$ we obtain
\begin{equation}\label{e:G(un-)<}
\int_\Omega|\nabla u_n^-|^2 + e^{u_n^-} + \int_\Gamma g|u_n^-|
 \leq 
 C|\Omega|.
\end{equation}
Note that it is easy to prove that the Poincar\'e inequality holds in 
$\left\{u-\frac{1}{|\Gamma|}\int_\Gamma u,\; u\in H^1(\Omega)\right\}$. Then from (\ref{e:G(un-)<}) it follows
that ${\displaystyle v_n=u_n^- - \frac{1}{|\Gamma|}\int_\Gamma u_n^-}$ is bounded in $H^1(\Omega)$.
Therefore, from (\ref{e:G(un-)<}) it follows $u_n^-=v_n+\frac{1}{|\Gamma|}\int_\Gamma u_n^-$ is bounded in $H^1(\Omega)$, which proves  that
$u_n$ is bounded in $H^1(\Omega)$.\\
ii)
{\it Up to a subsequence, $u_n$ converges to $u\in H^1(\Omega)$ weakly in $H^1(\Omega)$,
strongly in $H^{1-s}(\Omega)$, $s\in(0,1]$ and almost everywhere in $\Omega$, $u$ solves (\ref{e:min-G}) and
\begin{eqnarray}
 \int_\Omega|\nabla u|^2 &\leq&\liminf_{n\to\infty}\int_\Omega|\nabla u_n|^2 \leq C,	\label{e:|un|H1}\\
 \int_\Gamma g u_n &=&\lim_{n\to\infty}\int_\Gamma g_n \leq C,				\label{e:gun}\\
 \int_\Omega u^2 &=&\lim_{n\to\infty}\int_\Omega|u_n|^2 \leq C,			\label{e:|un|L2}\\
 \int_\Omega e^u &\leq& \liminf_{n\to\infty}\int_\Omega e^{u_n} \leq G(0)=|\Omega|,		\label{e:|eu|L1}
\end{eqnarray}
with $C=C(N,\Omega)$}.
Indeed, up to a subsequence, $u_n$ converges to a certain $u\in H^1(\Omega)$ weakly in $H^1(\Omega)$,
strongly in $H^{1-s}(\Omega)$, $s\in(0,1]$ and almost everywhere in $\Omega$. Then (\ref{e:|un|H1})-(\ref{e:|eu|L1})
follows. Note that the estimate (\ref{e:|eu|L1}) follows from Fatou's lemma because 
$e^{u_n}\geq0$, $\ds{\sup\left\{\int_\Omega e^{u_n},\ n\in\mathbb N\right\}<\infty}$ and $e^{u_n}$ converges to $e^{u}$ almost everywhere in $\Omega$.
Therefore, $u\in H^1(\Omega)$ and solves (\ref{e:min-G}).\\
iii)
{\it For any $\ds{v\in H^1(\Omega)\cap\left(\{v\in L^\infty(\Omega) \cup \{v\geq0\}\cup \{u\}\right)}$ we have}
\begin{equation}\label{e:G'(v)=0}
G'(u)(v)= \int_\Omega\nabla u\cdot\nabla v + e^u v - \int_\Gamma g v=0.
\end{equation}
Indeed, let first $v\in H^1(\Omega)\cap L^\infty(\Omega)$. 
Note that from (\ref{e:|eu|L1}), $e^u\in L^1(\Omega)$ and then $e^{u+tv}=e^u e^{tv}\in L^1(\Omega)$, for 
all $t\in\mR$. From $G(u+tv)-G(u)\geq0$ we obtain
\[ 
\sgn(t)
\frac{G(u+tv)-G(u)}{t}
=
\sgn(t)
\left(
\int_\Omega 
\nabla u\cdot\nabla v + \frac{t}{2}|\nabla v|^2
+
\frac{e^{tv}-1}{t} e^u
-
\int_\Gamma 
v
\right)
\geq0,
\]
which, after passing in the limit as $t$ tends to zero, yields (\ref{e:G'(v)=0}).

Taking in  (\ref{e:G'(v)=0}) $v=u_n:={\rm sign}(u)\min\{|u|,n\}$, $n\in\mathbb N$, gives
\begin{equation}\label{e:G'(u)(un)}
\int_\Omega|\nabla u_n|^2 + (e^u -1) u_n = \int_\Gamma g u_n - \int_\Omega u_n.
\end{equation}
We can pass in the limit in (\ref{e:G'(u)(un)}) as $n$ tends to infinity. Note that  
as $(e^u-1)u_n\geq0$, (\ref{e:G'(u)(un)}) implies  $\ds{\sup\left\{\int_\Omega (e^u-1)u_n,\ n\in\mathbb N\right\}<\infty}$.
As $(e^u-1)u_n$ is increasing, from the monotone convergence theorem of 
Beppo Levi, it follows ${\displaystyle\lim_{n\to\infty}\int_\Omega (e^u-1)u_n=\int_\Omega (e^u-1)u}$. 
Therefore, (\ref{e:G'(u)(un)}) yields
\begin{equation}\label{e:G'(u)=0}
 \int_\Omega |\nabla u|^2 + e^u u - \int_\Gamma g u 
 =
 0.
\end{equation}
Finally, let $v\in H^1(\Omega)\cap\{v\geq0\}$.
Then (\ref{e:G'(v)=0}) holds for $v_n=\min\{v,n\}$ instead of $v$.
%, so 
%\begin{equation}\label{e:G'(v)=0++}
%\int_\Omega\nabla u\cdot\nabla v_n + e^u v_n = \int_\Gamma g v_n.
%\end{equation}
Using the same argument as for the case $v=u$ and passing in the  limit as $n\to\infty$ in the equality  above 
yields (\ref{e:G'(v)=0}).\\
iv)
{\it The solution $u$ is bounded from above and}
\begin{equation}\label{e:G'(v)=0*}
 \int_\Omega \nabla u\cdot\nabla v + e^u v = \int_\Gamma g v,\quad
 \forall v\in H^1(\Omega).
\end{equation}
To prove this, we follow the technique due to {Stampacchia} \cite{stampacchia-1} as follows.
For $k\in\mathbb N$ set $v_k=\max\{u-k,0\}$, $A_k=\{x\in\Omega,\; v_k>0\}=\{x\in\Omega,\ u>k\}$. 
From iii) we can take $v=v_k$ in (\ref{e:G'(v)=0}), so
\[
 \int_\Omega|\nabla v_k|^2 + e^u v_k = \int_\Gamma g v_k
 \leq
 |g|\|v_k\|_{W^{1,1}(A_k)}
 \leq
 C |g||A_k|^{1/2}\|v_k\|_{H^1(\Omega)},\qquad
 |A_k|=measure(A_k),
\]
which, as $e^u v_k\geq v_k^2$, implies
\begin{eqnarray}
\int_\Omega|\nabla v_k|^2 + |v_k|^2 
&\leq&
C
|A_k|,\quad
C=C(g).							\label{e:|vk|H1}
\end{eqnarray}
We recall the Sobolev embedding 
$H^1(\Omega)\subset L^{2^*}(\Omega)$, where $\frac{1}{2^*}=\frac{1}{2}-\frac{1}{N}$ for $N\geq 3$ and 
$H^1(\Omega)\subset L^q(\Omega)$, $q\in[1,\infty)$ for $N=2$.
Now let $h>k$. 
Then from (\ref{e:|vk|H1}), we obtain
\begin{eqnarray}
(h-k)^2|A_h|^{\frac{2}{2^*}}
&\leq&
\left(\left(\int_{A_k}|v_k|^{2^*}\right)^{\frac{1}{2^*}}\right)^2
\leq
C\|v_k\|_{H^1(\Omega)}^2
\leq
C|A_k|,\quad \mbox{\it and so}					\label{e:(h-k)^2}\\
|A_h|
&\leq&
\frac{C}{(h-k)^{2^*}}|A_k|^{\frac{2^*}{2}},\qquad
\frac{2^*}{2}=\frac{N}{N-2}>1,
\quad
\mbox{\it if $N\geq 3$},
\end{eqnarray}
and
\begin{eqnarray}
(h-k)^2|A_h|^{\frac{2}{q}}
&\leq&
\left(\left(\int_{A_k}|v_k|^{q}\right)^{\frac{1}{q}}\right)^2
\leq
C\|v_k\|_{H^1(\Omega)}^2
\leq
C|A_k|,\quad \mbox{\it and so}					\label{e:(h-k)^2,N=2}\\
|A_h|
&\leq&
\frac{C}{(h-k)^{q}}|A_k|^{\frac{q}{2}},\qquad
\mbox{\it for all\;\, $q\geq1$},
\quad
\mbox{\it if $N=2$}.
\end{eqnarray}
The conditions of Lemma B.1, \cite{stampacchia-1}, are fulfilled (with $q>2$ in $N=2$)
($\varphi(h)\leq (C/(h-k)^\alpha)\varphi(k)^\beta$, $C, \alpha>0$, $\beta>1$, $h>k\geq0$, $\varphi(h)=|A_h|$). 
Then $|A_h|=0$ in $[h_0,\infty)$, for a certain $h_0$, which proves that $u$ is bounded from above.

As a corollary, (\ref{e:G'(v)=0*}) follows from (\ref{e:G'(v)=0}).\\
v)
{\it We have $u\in W^{2,q}(\Omega)$, $q\in(1,\infty)$}.
Indeed, note that $u\in H^1(\Omega)$ is a weak solution of 
\[
 -\Delta u = f:=-e^u\;\; in\;\, \Omega,\quad
 \partial_\nu u = g\chi_\Gamma\;\, on\; \partial\Omega,\quad
  f\in L^\infty(\Omega)\subset L^q(\Omega),\;\, q\in[1,\infty].
\]
We recall the regularity results for $-\Delta$ with Neumann boundary conditions, see \cite{adn-1+2}.
As our domain $\Omega$ is not $C^2$, we consider $\hat{\Omega}$, the domain obtained by reflecting $\Omega$ with respect to the plane $\{x_1=0\}$, and then by extending it $2\ell$ periodically along the $x_1$ axis.
Let $\hat{\Gamma}$ be its boundary,  which is $C^2$ owing to the assumption on $\Gamma_0$ and $\Lambda$.
Finally, let $\hat{u}$ be the extension of $u$ in $\hat{\Omega}$ obtained by reflection with respect to $\{x_1=0\}$ and then by extending it $2\ell$ periodically along the $x_1$ direction.
It is easy to show that
\begin{equation}
 -\Delta\hat{u} = \hat{f}:=-e^{\hat{u}}\;\; in\;\, \hat{\Omega},\quad
  \partial_{\hat{\nu}}\hat{u} = g\;\, on\; \partial\hat{\Omega},\quad
  \hat{f}\in L^\infty(\hat{\Omega})\subset L^q_{loc}(\hat{\Omega}),\;\, q\in[1,\infty],
\end{equation}
where $\hat{\nu}$ is the unit normal vector to $\partial\hat{\Omega}$ outward to $\hat{\Omega}$.
From $W^{2,q}$ regularity of $-\Delta$ we have $u\in W^{2,q}(\Omega)$, $q\in(1,\infty)$ 
(see \cite{adn-1+2, brezis-1}) and 
\begin{equation}\label{e:|u|W2q<}
\|u\|_{W^{2,q}(\Omega)}\leq C(\Omega)(\|e^u\|_{L^q(\Omega)} + \|u\|_{L^q(\Omega)}+\|g\|_{W^{1-1/q,q}(\Gamma)}).
\end{equation}
As a corollary, as $q>1$ is arbitrary, we get $u\in C^0(\ol{\Omega})$.\\
vi)
Finally, for the uniqueness, we point out that if $u$ and $\hat{u}$ are two solutions of (\ref{e:uw}), then
\[
 \int_\Omega (\nabla (u-\hat{u})\cdot\nabla v) + u_0(e^{u}-e^{\hat{u}})v = 0.
\]
Taking $v=u-\hat{u}$ and using the monotonicity of $e^u$ gives $u-\hat{u}=0$.
\hfill$\Box$

\begin{remark}\label{r:|u|<infty}
From the proof of Theorem \ref{th:u-exist+uniq} (step iv)), we see that $u$ is bounded from above regardless of the regularity of $\Omega$.
The boundedness of $u$ (from below) is derived by using the regularity of $-\Delta$ (and the Sobolev embedding theorem),
which uses strongly the regularity of $\Omega$ 
\colre{(the assumption $\Gamma_0$ is $C^2$, (\ref{e:gamma0}) and (\ref{e:L}))}.
\end{remark}

\subsubsection{Shape differentiation of the state variables}\label{sss:shd-u}
Now we turn our attention to the differentiability with respect to $\lambda$ of \colre{$U^0=\lambda+V^0$, given by Theorem \ref{th:U}, and of $u$, given by theorem \ref{th:u-exist+uniq}}.

Note that for $\lambda\in\Lambda$ fixed such that $I+\lambda$ is invertible, $I+\mu\in\Lambda$ 
is invertible for $\mu\in\Lambda$ near $\lambda$ and $(I+\mu)(D_0)=D_0$.

We consider \colre{first the differentiability of the} function $\mu\in\Lambda \mapsto U^0(\mu)\in H^1(S_0,\mR^N)$ near $\mu=\lambda$. \colre{We have}
\begin{theorem}\label{th:U'}
Let $\lambda\in \Lambda$ be given and $U^0(\mu)$ be the weak solution of (\ref{e:U0}), (\ref{bc:U0-Z0})-(\ref{bc:U0-G0}) as in Theorem \ref{th:U}, for $\mu$ near $\lambda$.
Then the map 
$\mu\mapsto U^0(\mu)$ is $C^1$ near $\lambda$ from 
$\Lambda$ to $H^1(S_0,\mR^N)$.

Furthermore, if $\dot{U}^0$ is the derivative of $U^0(\mu)$ at $\lambda$ in the direction $\mu$, then 
it satisfies
\begin{eqnarray}
- \nabla \cdot \epsilon^0(\dot{U}^0)
&=&
0\quad in\ S_0, 			\label{e:tU0}	\\
\dot{U}^0
&=&
0\quad on\;\;  I_0\cup\Sigma_0\cup\Pi_0,		\label{bc:tU0-I0}\\
\dot{U}^0
&=&
\mu\quad on\ \Gamma_0.			\label{bc:tU0-G0}
\end{eqnarray}
\end{theorem}
{\bf Proof}.
\colre{The proof of this theorem is straightforward and we will not present it here}.
\hfill$\Box$

\begin{theorem}\label{th:u'}
\colre{Assume (\ref{e:gamma0}) and let} $\lambda\in \Lambda$ with $I+\lambda$ invertible and $u(\mu)\in W^{2,q}(\Omega_\mu)$ be the solution of (\ref{e:uw}), as given by Theorem \ref{th:u-exist+uniq}, for $\mu\in\Lambda$ near $\lambda$.
Then the function $\mu\mapsto u_\mu:=u(\mu)\circ(I+(\mu-\lambda)\circ(I+\lambda)^{-1})$ is differentiable near
$\lambda$ from $\Lambda$ to $W^{2,q}(\Omega_\lambda)$, \colre{$q>1$}.

If $u'$ is the shape derivative of $u(\mu)$ at $\lambda$ in the direction $\mu\in\Lambda$ \colre{(see \cite{simon-1})}, we have
\begin{eqnarray}
u_\mu'(\lambda)(\mu)
&=&
u'+\mu\cdot\nabla u(\lambda)\quad in\;\; W^{1,q}(\Omega_\lambda),			\label{e:hatu'}\\
-\Delta u'+ \colre{u_0}e^{u(\lambda)} u'
 &=&
 0\;\; in\;\, \Omega_\lambda.					\label{e:Deltau'}
\end{eqnarray}
\colre{Furthermore, if $u(\lambda)\in W^{3,q}(\Omega_\lambda)$ then
$\partial_{\nu_\lambda} u'\in W^{1-1/q,q}(\partial\Omega_\lambda)$ and }
\begin{eqnarray}
 \partial_{\nu_\lambda} u'
 &=&
 0,\;\; on\;\, I_0\cup O_0,					\label{e:Dnu'-I0O0}\\
 \partial_{\nu_\lambda} u'
 &=&
 (\partial_{\tau_\lambda} u(\lambda)\cdot\partial_{\tau_\lambda}(\mu\cdot\nu_\lambda) 
 - 
 \partial^2_{\nu_\lambda}u(\lambda)(\mu\cdot\nu_\lambda))\chi_{\Gamma_\lambda}\;\; on\;\, \colre{\Gamma_\lambda}.		\label{e:Dnu'}
\end{eqnarray}
Here $\partial_{\tau_\lambda}(\cdot)=\nabla(\cdot)-(\nu_\lambda\cdot\nabla(\cdot))\nu_\lambda$ is the tangential gradient on $\Gamma_\lambda$.
\end{theorem}
{\bf Proof}.
First, we prove that $\mu\mapsto u_\mu$ is differentiable from $\Lambda$ to $W^{2,q}(\Omega_\lambda)$ near $\lambda$.
For this, we follow the classical approach  of proving the shape differentiability of boundary value problems,
which is based on the implicit function theorem, see for example \cite{simon-1}. 

As the domain $\Omega_{\mu}$ is not $C^2$, 
we consider an extension $\hat{\Omega}_{\mu}$ of $\Omega_{\mu}$, $\hat{\Gamma}_{\mu}$ its boundary and
the spaces 
$W^{2,q}_{2\ell}(\hat{\Omega}_\mu)$, 
$L^q_{2\ell}(\hat{\Omega}_\mu)$,
$W^{1-1/q,q}_{2\ell}(\hat{\Gamma}_\mu)$,
and
$\hat{u}(\mu)\in W^{2,q}_{2\ell}(\hat{\Omega}_\mu)$  the extension of $u(\mu)$ in $\hat{\Omega}_{\mu}$ 
\colre{(as in Theorem \ref{th:u-exist+uniq})}.

It is easy to point out that  
\[
 -\Delta \hat{u}(\mu) \colre{+u_0e^{\hat{u}(\mu)}=0}\;\; in\;\, \hat{\Omega}_{\mu},\quad
 \partial_{\hat{\nu}_{\mu}}\hat{u}(\mu)=g\;\; on\;\, \hat{\Gamma}_{\mu},
\]
where $\hat{\nu}_{\mu}$ is the  unit normal vector to $\hat{\Gamma}_{\mu}$ exterior to
$\hat{\Omega}_{\mu}$.

It follows that 
$\hat{u}_\mu:=\hat{u}(\mu)\circ(I+\theta)\in W^{2,q}_{2\ell}(\hat{\Omega}_\lambda)$, 
\colre{where $\theta=(\hat{\mu}-\hat{\lambda})\circ(I+\hat{\lambda})^{-1}$},
and it satisfies (see \cite{simon-2})
\begin{eqnarray*}
 -\sum_{i,j=1,N} Q_{ij}(\mu)\partial_j(Q_{ik}(\mu)\partial_k \hat{u}_\mu) 
 + 
 u_0 e^{\hat{u}_\mu}
 &=&
 0\;\; in\;\; \hat{\Omega}_\lambda,	\\
 -\frac{1}{|Q\cdot\nabla d_\lambda|}
 \nabla d_\lambda\cdot {^tQ}\cdot Q\cdot\nabla \hat{u}_\mu &=&g\;\; on\;\; \hat{\Gamma}_\lambda,
\end{eqnarray*}
where 
\colre{$Q_{ij}(\mu)={^t[\nabla(I+\theta)]^{-1}}$},
$d_\lambda(x)=(1-2\chi_{\hat{\Omega}_\lambda}(x))dist(x,\hat{\Gamma}_\lambda)$ 
\colre{is the oriented distance to $\hat{\Gamma}_{\lambda}$},
and \colre{$\hat{\lambda}$ and $\hat{\mu}$} are the extensions  in $\hat{\Omega}_0$ of \colre{$\lambda$ and $\mu$ respectively} as in Theorem \ref{th:u-exist+uniq}.

Then we use (as it is standard) the implicit function theorem. Namely, we consider
\begin{eqnarray*}
T=(A,B)&:&\Lambda\times W^{2,q}_{2\ell}(\hat{\Omega}_\lambda) 
\mapsto 
L^q_{2\ell}(\hat{\Omega}_\lambda)\times W^{1-1/q,q}_{2\ell}(\hat{\Gamma}_\lambda),\\
A(\mu,\hat{v})
&=&
 -\sum_{i,j=1,N} Q_{ij}(\mu)\partial_j(Q_{ik}(\mu)\partial_k \hat{v})  +  u_0e^{\hat{v}},\\
B(\mu,\hat{v})
 &=&
 \frac{1}{|Q\cdot\nabla d_\lambda|} (\nabla d_\lambda\cdot{^tQ}\cdot Q\cdot\nabla \hat{v} + g.
\end{eqnarray*}
Note that $T(0,\hat{u}_\lambda)=0$, and from \cite{simon-2}, $T$ is $C^1$ near $0$. Furthermore,
\begin{equation}\label{e:DuF}
\partial_{\hat{v}} T(0,\hat{u}_\lambda)(\hat{v})
=
(-\Delta \hat{v}+u_0\colre{\hat{v}}e^{\hat{u}_\lambda} ,\partial_{\hat{\nu}_\lambda}\hat{v})
\in
L^q_{2\ell}(\hat{\Omega}_\lambda)\times W^{1-1/q,q}_{2\ell}(\hat{\Gamma}_\lambda).
\end{equation}
\colre{But from the $W^{2,q}$ regularity of $-\Delta+I$, see \cite{adn-1+2, brezis-1}},
$\partial_v T(0,\hat{u}_\lambda)$ defines an isomorphism from $W^{2,q}_{2\ell}(\hat{\Omega}_\lambda)$ to 
$L^q_{2\ell}(\hat{\Omega}_\lambda)\times W^{1-1/q,q}_{2\ell}(\hat{\Gamma}_\lambda)$.
Then, the differentiability of 
$\mu\mapsto \hat{u}_\mu$, and so of $\mu\mapsto u_\mu$, follows by using the implicit function theorem.

Furthermore, (\ref{e:hatu'}), (\ref{e:Deltau'}), (\ref{e:Dnu'}) follow from \cite{simon-1} and \cite{simon-2}.
\hfill$\Box$

\subsubsection{Shape differentiation of the energy (Proof of Proposition \ref{p:E'=})}\label{sss:shd-E}
From Theorem \ref{th:U'}, it follows that $\mu\mapsto E_{mech,s}(\mu)$ is differentiable at $\lambda$,
see \cite{simon-1}. Then
%%If $E'_{mech,s}(\lambda)(\mu)$ denotes the derivative of $\mu\mapsto E_{mech,s}(\mu)$ at $\lambda$ in the direction %%$\mu$, we have
\begin{eqnarray*}
 \colre{\frac{d}{d\mu}}E_{mech,s}(\lambda)(\mu)
 &=&
 \frac{1}{2}\int_{S_0}
 \epsilon^0(\dot{U}^0):[\nabla U^0] + \epsilon^0(U^0):[\nabla \dot{U}^0]	
 -
\int_{S_0} p_S^0\nabla\cdot \dot{U}^0,
\end{eqnarray*}
where $\dot{U}^0$ is given by (\ref{e:tU0})-(\ref{bc:tU0-G0}).

Note that for $U, V\in H^1(S_0,\mR^N)$ we have 
$
 {\displaystyle
 \int_{S_0}
\epsilon^0(U):[\nabla V] 
= 
\int_{S_0}
\epsilon^0(V):[\nabla U] 
}$.
Then from (\ref{e:U0}), (\ref{bc:U0-Z0})-(\ref{bc:U0-G0}), (\ref{e:tU0})-(\ref{bc:tU0-G0}), and  (\ref{e:sigma.n=sigma0.n0}), we obtain
\begin{eqnarray}
\colre{\frac{d}{d\mu}}E_{mech,s}(\lambda)(\mu)
&=& 
\int_{S_0}
 \epsilon^0(U^0):[\nabla \dot{U}^0]
 -
 p_S^0\nabla\cdot \dot{U}^0		\nonumber\\
&=&
-\int_{S_0}(\nabla\cdot\sigma^0(U^0))\cdot \dot{U}^0
+
\int_{\Gamma_0} \left(\sigma^0(U^0)\cdot\nu_0^S\right)\cdot \dot{U}^0d\Gamma_0	\nonumber\\
&=&
-
\int_{\Gamma_\lambda} 
\left(\sigma(U^0)\cdot\nu_\lambda\right)\cdot\mu d\Gamma_\lambda,		\label{e:Emech,S'}
\end{eqnarray}
where \colre{$\nu_0^S$} is the unit normal vector to $\partial S_0$ exterior to $S_0$.

From $p=RT(c-c_0)+p_0$ we get
\begin{equation}\label{e:Emech,O+Eel}
\hspace*{-7mm}
E_{mech,l}(\mu)+E_{el}(\mu)
=
-
\int_{\Omega_{\mu}}
\left(
\frac{\varepsilon}{2}|\nabla \varphi(\mu)|^2 + RTc(\mu)
\right)
-
\int_{\Omega_{\mu}}(p_0-RTc_0)
-
\int_{\Gamma_{\mu}}\sigma_c\varphi,
\end{equation}
where 
$\ds{\varphi=\varphi(\mu)=-\frac{RT}{F}u(\mu)}$,  
$\ds{c=c(\mu)=e^{-\frac{F}{RT}\varphi(\mu)}}$.
Note that from Theorem \ref{th:u'}, the map 
$\mu\in\Lambda \mapsto c_\mu:=c(\mu)\circ(I+(\mu-\lambda)\circ(I+\lambda)^{-1})\in W^{2,q}(\Omega_\lambda)$ is $C^1$ 
near $\lambda$.
Furthermore, if $\varphi'$, resp. $c'$, is the shape derivative at $\lambda$ in the direction $\mu$ 
of $\varphi$, resp. $c$, then
$\ds{\varphi'=-\frac{RT}{F}u'}$, $\ds{c'=-\frac{F}{RT}c\varphi'=cu'}$, see  \cite{simon-1}.
From Theorem \ref{th:u'} and \cite{simon-1}, \cite{simon-2} we obtain
\begin{eqnarray}
\frac{d}{d\mu}
\left(\int_{\Omega_{\mu}}
\frac{\varepsilon}{2}|\nabla \varphi(\lambda+\mu)|^2 + RTc\right)(\lambda)(\mu)
&=&
\int_{\Omega_{\lambda}}
(\varepsilon\nabla \varphi\cdot\nabla\varphi' -  Fc\varphi') 	\nonumber\\
&+&
\left(
\int_{\Gamma_\lambda}\frac{\varepsilon}{2}|\nabla\varphi(\lambda)|^2
+
RTc\right)
(\mu\cdot\nu_\lambda),							\label{e:Emech,O'}\\
\frac{d}{d\mu}
\left(\int_{\Omega_{\mu}}(p_0-RTc_0)\right)(\lambda)(\mu)
&=&
\int_{\Gamma_\lambda}(p_0-RTc_0)(\mu\cdot\nu),				\label{e:(p0-RTc0)'}\\
\frac{d}{d\mu}
\left(\int_{\Gamma_{\mu}}\sigma_c\varphi\right)(\lambda)(\mu)
&=&
\int_{\Gamma_\lambda}
\sigma_c\varphi'
+
\sigma_c(\mathcal H_\lambda\varphi + \partial_{\nu_\lambda}\varphi)
(\mu\cdot\nu),								\label{e:intGphi'}
\end{eqnarray}
where  ${\cal H}_\lambda$ is the mean curvature of $\Gamma_\lambda$. 

%%Let $E'_{mech,l}(\lambda)(\mu)$, resp. $E_{el}'(\lambda)(\mu)$, be the shape derivative of 
%%$\mu\mapsto E_{mech,l}(\mu)$, resp. $\mu\mapsto E_{el}(\mu)$, at $\lambda$ in the direction $\mu\in\Lambda$. 

From (\ref{e:Emech,O+Eel}), (\ref{e:Emech,O'}), (\ref{e:(p0-RTc0)'}), (\ref{e:intGphi'}) we get
\begin{eqnarray}
\colre{\frac{d}{d\mu}}E_{mech,l}(\lambda)(\mu)
+
\colre{\frac{d}{d\mu}}E_{el}(\lambda)(\mu)
&=&
\int_{\Omega_\lambda}
-
\varepsilon(\nabla\varphi\cdot\nabla\varphi')
+
Fc\varphi' 
-
\int_{\Gamma_\lambda}\sigma_c\varphi'			\nonumber\\
&-&
\int_{\Gamma_\lambda}
\left(
\frac{\varepsilon}{2}|\nabla\varphi|^2
+ 
(RT(c-c_0)+p_0) 
+
\sigma_c(\mathcal H\varphi + \partial_\nu\varphi)
\right)
(\mu\cdot\nu)						\nonumber\\
&=&
\int_{\Omega_\lambda}
(\nabla\cdot(\varepsilon\nabla\varphi)+Fc)\varphi'
+
\int_{\partial\Omega_\lambda}(-\varepsilon\partial_{\nu_\lambda}\varphi-\sigma_c\chi_\Gamma) \nonumber\\
&-&
\int_{\Gamma_\lambda}
\left(
\frac{\varepsilon}{2}|\nabla\varphi|^2
+ 
p 
+
\sigma_c \mathcal H\varphi -\varepsilon(\partial_\nu\varphi)^2
\right)
(\mu\cdot\nu)						\nonumber\\
&=&
-
\int_{\Gamma_\lambda}
\left(
p
+
\frac{\varepsilon}{2}
\left(|\partial_\tau\varphi|^2-|\partial_\nu\varphi|^2\right)
+
\sigma_c \mathcal H\varphi
\right)
(\mu\cdot\nu),					\label{e:Emech,O'+Eel'}
\end{eqnarray}
because 
$\ds{\frac{1}{\varepsilon}\sigma^2_c=\varepsilon|\partial_\nu\varphi|^2}$,
where $\partial_\tau\varphi$ is the tangential gradient of $\varphi$ on $\Gamma_\lambda$.

For $E_{st}$, 
%%if $E'_{st}(\lambda)(\mu)$ is the derivative of $\mu\mapsto E_{st}(\mu)$ at $\lambda$ in the direction $\mu\in\Lambda$, 
from classical shape calculus we have
\begin{equation}\label{e:Est'}
 \colre{\frac{d}{d\mu}}E_{st}(\lambda)(\mu)=\int_\Gamma \gamma \mathcal H (\mu\cdot\nu).
\end{equation}

Adding (\ref{e:Emech,S'}), (\ref{e:Emech,O'+Eel'}), (\ref{e:Est'}) gives (\ref{e:E'(l)(m)}).

Finally, (\ref{e:fb-v}) follows from $\colre{\frac{d}{d\mu}}E(\lambda)(\mu)=0$ and the arbitrariness of $\mu$.%
\hfill$\Box$

\section{Numerical results}\label{s:nr}
In this section, we will present approximations of the equilibrium interface $\Gamma$, a solution of the free boundary interface (\ref{e:fb-v+}), resp. (\ref{e:fb}), \colre{in 2D and 3D}, by using the algorithms (\ref{e:gm}), resp. (\ref{e:fpm}).

Our algorithms have been implemented in the commercial finite element software, COMSOL 3.4.
Since the Young-Laplace equation ((\ref{e:fb-v+}) and (\ref{e:fb}))
%our shape gradient  (\ref{e:E'(l)(m)})
 depends on the mean curvature $\mathcal H$ of the interface, 
 a critical aspect regarding a good approximation of the interface $\Gamma$ 
is an accurate computation of $\mathcal H$, see \cite{MGCR}.
%Following the ideas in \cite{MGCR},
At first we calculate %in weak formulation
an extension $\mathcal V$ of the normal vector $\nu$ by
\begin{eqnarray}
\Delta \mathcal V 	
&=& 
0\qquad\,\,\,\,  in \,\,  \Omega_\lambda, \\
{\mathcal V} 		
&=&
\nu\chi_{\Gamma_\lambda}	\quad on \,\, \partial\Omega_\lambda.	\label{e:nu_bc_out}
\end{eqnarray}
%Note that the homogeneous b.c.~(\ref{e:nu_bc_out}) on $\partial\Omega_\lambda\backslash\Gamma_\lambda$ guarantees that %the %periodic
%extension of $\mathcal V$ in $\hat{\Omega}_\lambda$ is smooth.
If $\Gamma_\lambda$ is $C^{2,\alpha}$, $\alpha\in(0,1)$, by using the extension technique as in Theorem \ref{th:u-exist+uniq}, we obtain ${\mathcal V}\in C^{1,\alpha}(\ol{\Omega}_\lambda)$, see \cite{gilbarg+trudinger-1}.
This allows to compute the mean curvature
$\mathcal H := \nabla_\tau \cdot \nu=\nabla\mathcal V - (\nabla\mathcal V\cdot\nu)\nu$ on $\Gamma_\lambda$.
As numerically $\mathcal H$ represents oscillations, instead of $\mathcal H$ we consider a smoothed mean curvature, still denoted by $\mathcal H$, solution of
\begin{eqnarray}
%- \varepsilon_s^2 \Delta\mathcal H %+ (\mathcal H - \nabla \cdot {\mathcal V})
%&=& 
%0 \quad in\;\, \Omega_\lambda,			\label{e:-Delta H}\\
- \varepsilon_s \Delta_\tau \mathcal H + \mathcal H 
&=&
\nabla_\tau \cdot {\mathcal V}
\quad on\;\, \Gamma_\lambda, %\partial\Omega_\lambda,
\qquad
\mathcal H=0\quad on\;\, \partial\Gamma_\lambda,
\end{eqnarray}
%\textcolor{blue}{(what bc for $\mathcal H$? I have added $\mathcal H=0$ as bc ...)}
where $\varepsilon_s$ is a suitable small smoothing parameter that depends on the number of mesh elements and $\Delta_\tau$ denotes the Laplace-Beltrami operator on $\Gamma_\lambda$.
%
%Note that (\ref{e:-Delta H}) is equivalent to 
%$-\varepsilon_s \Delta(\mathcal H - \nabla \cdot {\mathcal V}) + (\mathcal H - \nabla \cdot {\mathcal V})=0$ in $\Omega$. Then $\mathcal H - \nabla \cdot {\mathcal V}\in C^0(\ol{\Omega}_\lambda)$, see e.g. \cite{brezis-3},
%so $\mathcal H \in C^0(\ol{\Omega}_\lambda)$.

For the update of the deformed geometry in the algorithms (\ref{e:fpm}) or (\ref{e:gm}), the ALE (Arbitrian Lagrangian Eulerian) module provided by the software is used.
Besides, in order to avoid inverted mesh elements and to guarantee a certain mesh quality, we remesh the domain before Step 1 of the algorithms (\ref{e:gm}) or (\ref{e:fpm}), if
necessary. As a stopping criterion, we consider $| \sup_{x \in \Gamma_n} |\lambda_{n}(x)| - \sup_{x \in \Gamma_{n-1}} |\lambda_{n-1}(x)| | \leq 10^{-3}$ 
%%\textcolor{red}{(shouldn't we instead consider the stop condition $\|g\|_X<10^{-3}$, where $g$ is the (free boundary) equation we solve and $X$ is a space, for example $L^\infty(\Gamma_n)$? To me this is more natural)}
for both algorithms.

We emphasize the large disparity of the physical parameters, see Table \ref{t:data}.
Therefore, in order to improve the conditioning of the numerical problem we have implemented our equations in non-dimensionalized form.
%NOT NEEDED: Within good approximation we may compute %(\ref{e:sigma=sigma0}) -- 
%(\ref{e:fp}) in linear elasticity, where terms of higher order then $\nabla F = \nabla U$ are neglected. Thus $F \approx \mathbf{1}$ and $\det(\nabla F) \approx 1$.
%
%For important parameters of our model see Table \ref{t:data}.
 References and remaining data needed for our simulation 
can be found in \cite{LBKN-11}. %The excess pressure $p_C := p_0 - p_S^0 + p^{ref}$ is an important parameter as we see in the following.
Note that, strictly speaking, $\sigma_c$ will change as the pore wall is deformed and sulfonic 
acid groups rearrange. However, it is beyond the scope of this contribution to model this change
also.

\begin{table}[!ht] %[htpb]
\centering
{
%\footnotesize
\begin{tabular}{%|
lccc%l
}
%\textbf
{Description} & %\textbf
{Symbol} & %\textbf
{Value} & %\textbf
{Unit} %& %\textbf
%{Reference}
 \\[3pt]
%\hline
Temperature & $T$ & $353$ & K \\
Surface tension of H$_2$O & $\gamma$ & $6.5 \times 10^{-2}$ & N/m \\ %& \cite{Kim} \\
%%%%%Elastic Young modulus & $E_{Y}$ & $275 \times 10^6$ & N/m${}^2$ \\ %& \cite{Plazanet%,Satterfield
%} \\
%%%%%Poisson ratio & $\nu_{P}$ & $0.491$ & $1$ \\%& \cite{Plazanet} \\
%
%
Elastic bulk modulus & $k_S$ & $5.09 \times 10^9$ & N/m${}^2$ \\ %& \\
Elastic shear modulus & $G_S$ & $9.22 \times 10^7$ & N/m${}^2$ \\%& \\
%
%
%%%%Electric bulk permittivity %of H$_2$O
%%%% & $\overline\varepsilon_r$ & $80$ & $1$ \\%& \cite{%KreuerEtAl,
%%%%Paul%,LadipoEtAl
%%%%} \\
Electric permittivity (near the interface) %(spatially averaged)  %of H$_2$O
 & $\varepsilon_r$ & $80$%$45$
 & $1$ \\%& \cite{%KreuerEtAl}
%%%%Electric permittivity %of H$_2$O
%%%% at wall & $\varepsilon_w$ & $8$ & $1$ \\%& \cite{%KreuerEtAl,
%%%%%%%%Paul%,PaulPaddison2,LadipoEtAl
%%%%%%%%} \\
Surface charge density %at wall
 & $-\sigma_c$ & $0.16022$ & C/m${}^2$ \\ %& \cite{%BergLadipo,
%Qiao} \\
%Typical length scale
%%%%%Diameter of a channel (Ref.~config.) & $d$ & $2 \times 10^{-9}$ & m \\ %& \cite{%LadipoEtAl,SchmidtRohrAndChen, %%%% or not in captions ?
%QiaoAluru2003} \\
%%%%%Length of 2 channels (Ref.~config.)%segment
%%%%% & $l = 10 d$ & $2 \times 10^{-8}$ & m \\ %%%%& \\ Length of membrane & $l_m$ & $10^{-4}$ & m \\ %& \cite{LadipoEtAl} \\%
Reference %Typical
 solid pressure  & $p_S^0 %\overline{p} %\approx p^* - \gamma \mathcal{H}^R
$ & $%- 6.1
-6.5 \times 10^7$ & N/m${}^2$ \\ %%%%& %e.g.~
%%%%\cite{Satterfield%,He
%%%%} \\
%Reference 
Typical reference liquid pressure & $p_0 %p^*
$ & $\approx
 1.4%4
 \times 10^6$ & N/m${}^2$ \\ %%%%& \cite{BergEikerling} \\
%%%%Typical pressure drop & $%p_0^{(1)} - p_0^{(2)} =
%%%%p_\Delta = p^* \, %40 \mbox{bar}
%%%% l/l_m$ & $8 \times 10^2$ & N/m${}^2$ & \cite{%BergLadipo,
%%%%He} \\
Typical reference concentration & $c_0 %\overline{c}
$ & $\approx 4.7790 \times 10^2$ & mol/m${}^3$ %%%% & \cite{BergLadipo} \\
%%%%%Typ.~internal el.~potential & $\overline{\phi}_0$ & $0.1$ & V & \cite{LadipoEtAl} \\
%%%%Typ.~mech.~displacement & $\overline{U}$ & $1.2611 \times 10^{-11}$ & m &
\end{tabular}
}
\caption[Database]{\label{t:data} Overview of significant parameters for our problem.}
\end{table}

\subsection{Simulations in 2D}
\begin{figure}[!ht]
\centering
\hspace{4mm}
\includegraphics[width=140mm,height=70mm]
{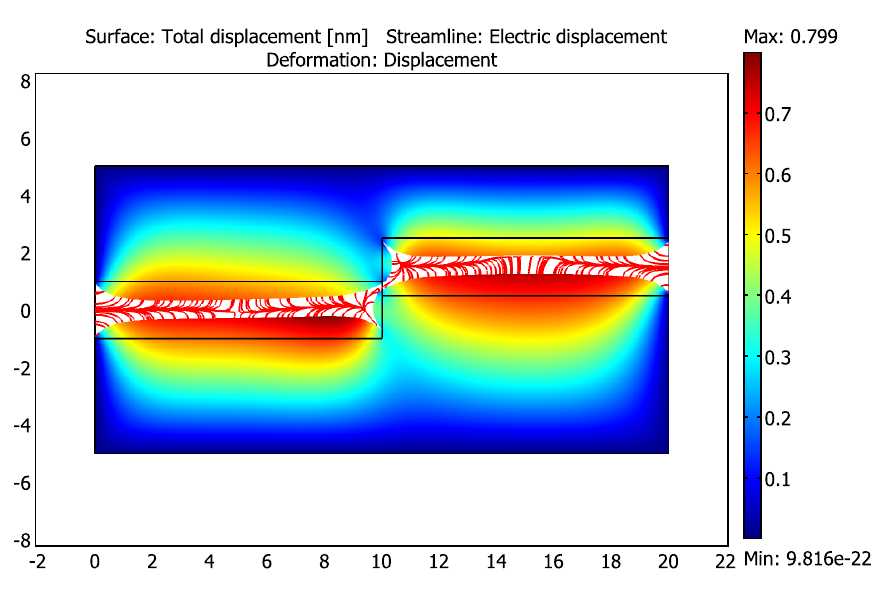}
\caption{\label{f:2d_closing} 
%$U_2$ component of the mechanical displacement field $U$
Euclidean norm of the %Absolute value of the
 mechanical displacement field, $| U^{(0)} |$, streamlines along the electric field $- \nabla \varphi$.
%($1$ unit $=1$ nm). 
%Negative 
Pressure $p_0 = 7348.96$ bar, modified Young-Laplace law. %(?)
A channel may close completely, the fixed point algorithm stops before the equilibrium pore shape has been determined. In the reference configuration, the cylinders have diameter $d = 2$ nm, length $l = 10$ nm and offset $s = 0.5$ nm.}
\end{figure}

\begin{figure}[!ht]
\centering
\begin{minipage}{150mm}
%\vspace*{2mm}
%\hspace*{4mm}
\includegraphics[width=70mm]
{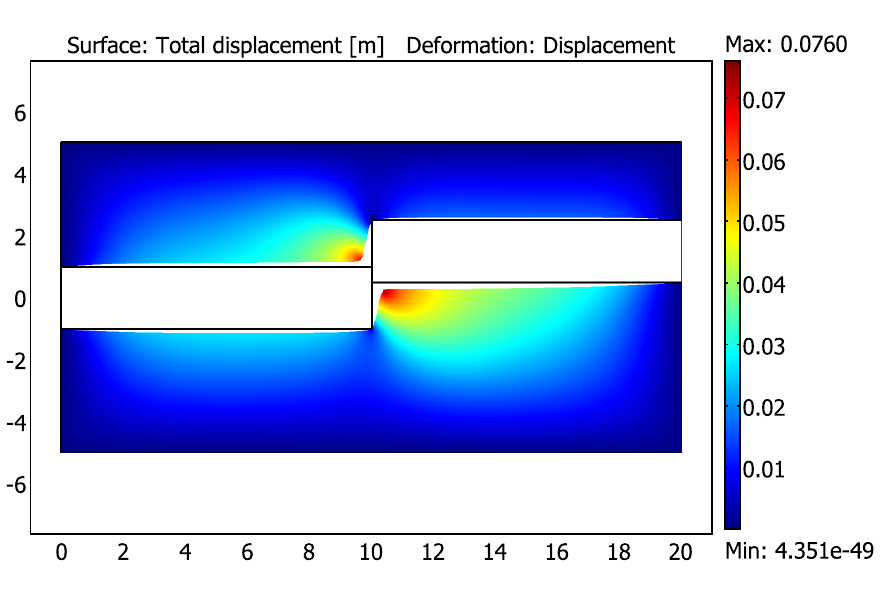}%\vspace*{-4mm}
\hspace*{2mm}
\includegraphics[width=70mm]
{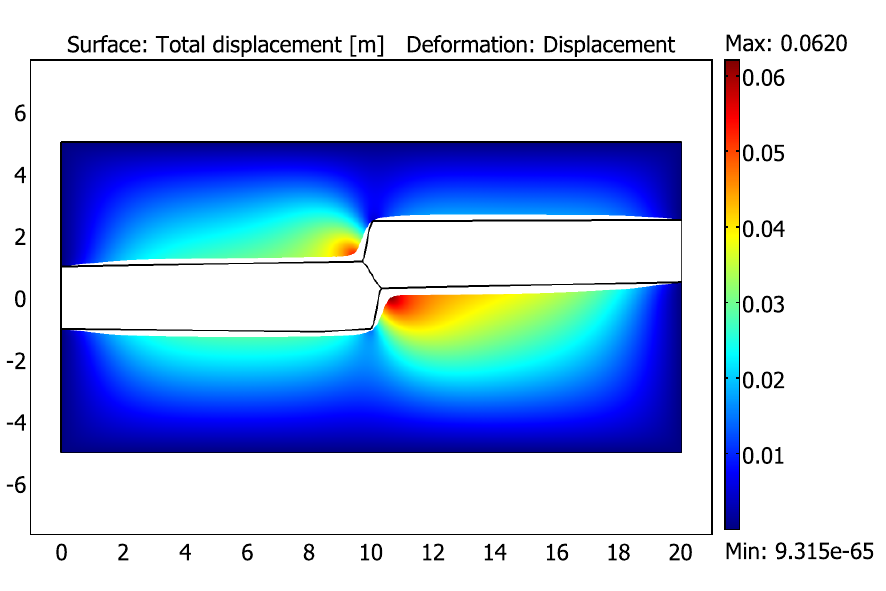}
\end{minipage}
\\%\vspace*{4mm}
\begin{minipage}{150mm}
%\vspace*{2mm}
%\hspace*{4mm}
\includegraphics[width=70mm]
{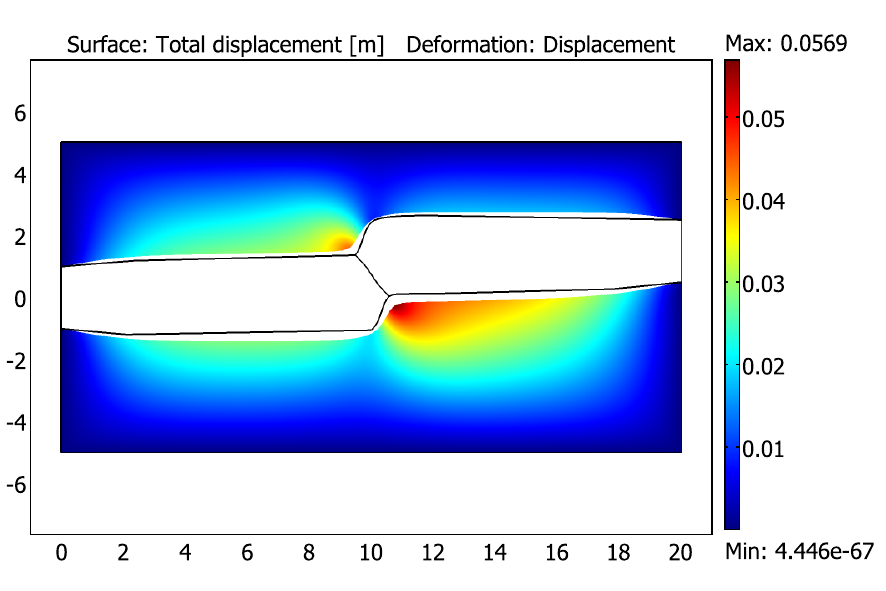}%\vspace*{-4mm}
\hspace*{2mm}
\includegraphics[width=70mm]
{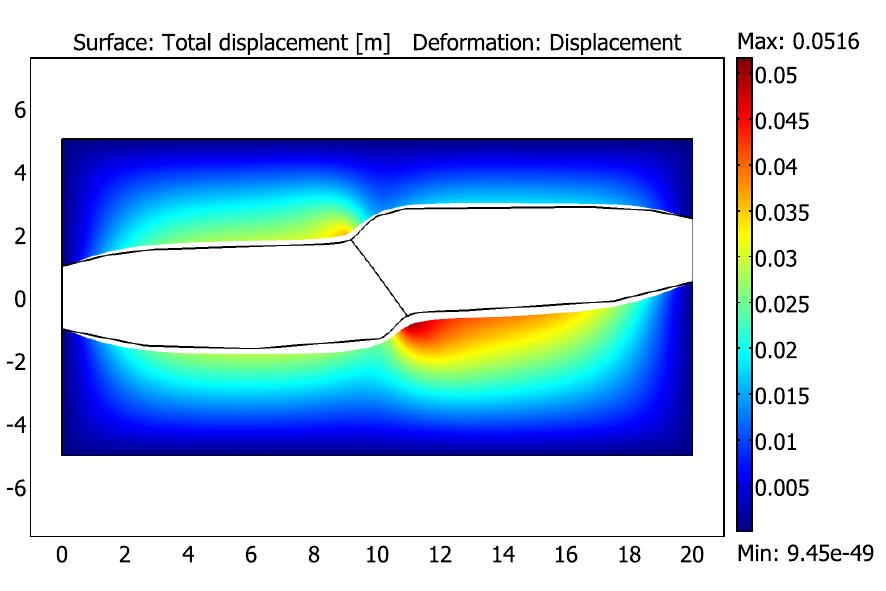}
\end{minipage}
\caption{\label{f:2d_conv} 
Euclidean norm of the mechanical displacement field, $|U^{(0)}|$, %($1$ unit $=1$ nm)
 at different steps. 
 Pressure $p_0 = 7349.03$ bar, classical Young-Laplace law. After $n = 5$ iterations (top left), $n = 10$ (top right), $n = 15$ (bottom left) and after numerical convergence, i.e.~$n = 28$ (bottom right). A channel expands until an equilibrium shape has been reached. Black lines indicate the shape of the initial or the deformed shape of the previous plot. In the reference configuration, the cylinders have diameter $d = 2$ nm, length $l = 10$ nm and offset $s = 0.5$ nm. %Here the ends of the channel at the r.h.s.~and at the l.h.s.~are fixed by modifying the mean curvature to $\mathcal H_{mod} = (\mathcal H - 1) \tanh(x_1) \tanh(20 - x) + 1$.
}
\end{figure}

Depending on $p_0$ (or $c_0$), the pore may close completely 
%without reaching an equilibrium shape before this event happens,
as in Fig.~\ref{f:2d_closing}, or expand  as in Fig.~\ref{f:2d_conv},
until an equilibrium shape has been achieved.

%\textcolor{blue}{
Fig. \ref{f:2d_closing} shows the shape of a closing pore, which is not a solution of (\ref{e:min(E)})
but rather a picture of the pore shape at the moment when the interface intersects itself (and when our algorithm is designed to stop). 
We think that the cusp-like corners might be due to electric repulsion, that prevents the                                               
two different parts of the membrane to join.
%}

In Fig.~\ref{f:2d_conv}, we show four plots at different stages of the convergence process for an expanding pore. 

Of course, the fixed point algorithm may be applied in the case of the modified Young-Laplace law, too. Then, the fixed point method (\ref{e:fpm}) and the variational method (\ref{e:gm}) yield similar results. Depending on an appropriate choice of the numerical parameter $k$, the variational method may converge faster, but it is very sensitive to perturbations.

\begin{figure}[!ht]
\centering
%\vspace*{2mm}
\hspace{4mm}
\includegraphics[width=140mm,height=70mm]{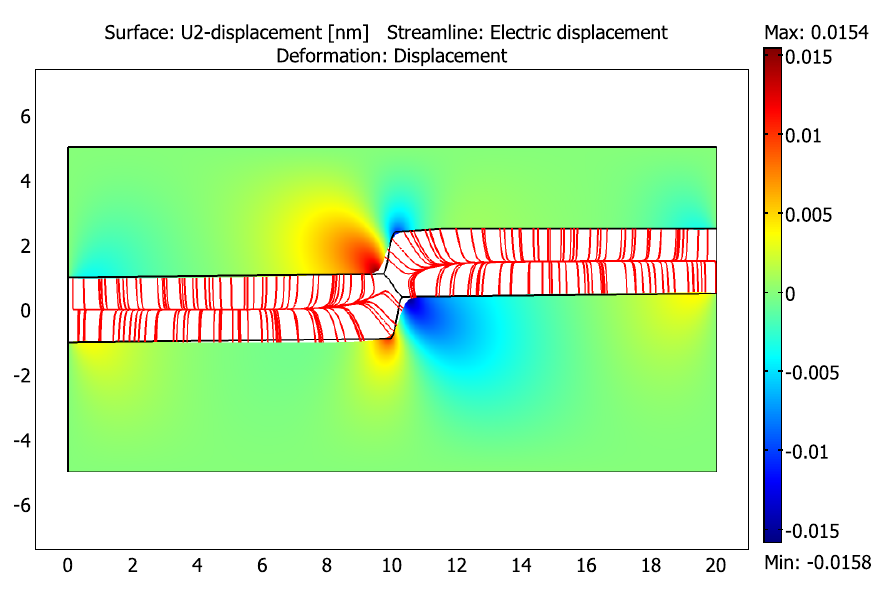}
\caption{\label{f:2d_mYL} 
$U_2^{(0)}$ component of the mechanical displacement field, %$U_2^{(0)}$,
 streamlines along the electric field $- \nabla \varphi$. $n = 10$ iterations, black lines show the shape after $n = 5$ iterations. Pressure $p_0 = 7349.35$ bar, modified Young-Laplace law. In the reference configuration, the cylinders have diameter $d = 2$ nm, length $l = 10$ nm and offset $s = 0.5$ nm.
}
\end{figure}

Results for the classical (Fig.~\ref{f:2d_conv}, top right) and the modified Young-Laplace law (see Fig.\ref{f:2d_mYL}), exhibit the significance of our modified Young-Laplace equation (\ref{e:fb-v}).

Note that for the two-dimensional case, the mesh consists typically of about $21,000$ elements.

%plot for closing

%4 plots for convergence, opening

%plot for modified YL

\subsection{Simulations in 3D}

Three dimensional simulations require many technical subtleties, e.g.~the initial mesh has to be chosen very fine near the interface but has to be sufficiently coarse otherwise so as to keep the number of variables small enough. For three space dimensions, the mesh consists typically of about $10^5$ %?
 elements.

\begin{figure}[!ht]
%\centering
%\begin{minipage}{80mm}%{150mm}
%\vspace*{2mm}
%\hspace*{4mm}
\centering
\hspace{4mm}
\includegraphics[width=140mm,height=70mm]{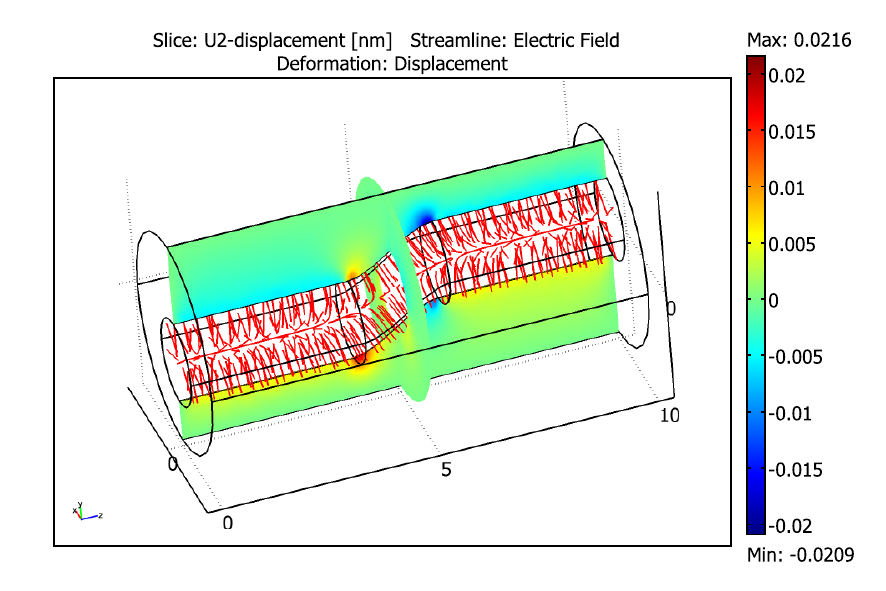}
%\vspace*{-4mm}
%\hspace*{2mm}
%%%%%%\includegraphics[draft,width=70mm]
%%%%%%{Figure5b}
%\end{minipage}
\caption{\label{f:3d_mYL}
$U_2^{(0)}$ component of the mechanical displacement field. Pressure $p_0 = 7349.03$ bar, %Left: 
modified Young-Laplace law. %Right: Modified Young-Laplace law in the presence of charges. In both plots
The cylinders have, in the reference configuration, diameter $d = 2$ nm, length $l = 5$ nm and offset $s = 0.5$ nm with a smoothed connection.}
\end{figure}

The situation of two channels with a smoothed connection in the reference configuration in case of the modified
 Young-Laplace law (variational method) is presented in Fig.~\ref{f:3d_mYL}.

3D simulations show that the interface solution of the modified Young-Laplace equation (obtained by using a gradient descent method) is slightly different  to the interface solution of the Young-Laplace equation (obtained with the fixed point method). 
%\textcolor{red}{``But slightly different equilibrium shapes are observed for the classical Young-Laplace law (by the %fixed point method)''.}
%\textcolor{blue}{shouldn't this prop. removed?}

%\textit{Open questions: Mesh independence, Speed of convergence}\\
%{\bf - we may consider these issues, but I think the main purpose of this paper is not on these issues; eventually we may add a short comment about them]}
%

From numerical experiments in 2D and 3D, we have seen that the convergence of the inner iteration loop (step 1.b) in both algorithms  is super-linear (by means of the SPOOLES solver), 
%\textcolor{blue}{(I still do not realize which is the 'inner iteration loop')}, 
with exception of the first iterations where a damped algorithm is essential. 
The outer loop, i.e. the iteration of Step 1 in order to update the interface shape, converges very slowly for 
the fixed-point method approach (see Fig.~\ref{f:2d_conv}).
Its convergence is faster with the variational method algorithm, but still slower than linear. 
In both cases, the convergence of the algorithm (meaning the outer loop) depends strongly on the initial offset $s$ between both pores.

\subsection{Discussion of a modified Young-Laplace law}
We have considered a model describing the charged fluid-elastomer interaction and the free boundary equation
(modified Young-Laplace equation) characterizing the fluid/elastomer equilibrium interface, see (\ref{e:fb-v+}). This equation differs from the classical Young-Laplace equation (\ref{e:fb}), commonly used in applications, see Remark \ref{r:extra-term}.

Our numerical computations show that the shape of the interface $\Gamma$, characterized by the modified Y-L equation 
(\ref{e:fb-v+}), is quite sensitive towards the reference values of pressures and concentration.
On the other side, the literature on parameter values for the surface tension between protonated water and elastomers is sparse.
The values of these parameters we have used in our computations, have been chosen such  that they lead to physically reasonable results. The difference between our parameters and the values from the literature might be explained by the modified Young-Laplace law (\ref{e:fb-v}) (see Remark \ref{r:extra-term}).
This law contains electrostatic terms, which are ignored in standard applications. 
We emphasize that generalizations for the surface tension %Young-Laplace law
 in the presence of surface charges can be found in the literature \cite{SchmicklerSantos}.

%Comparison of parameters
Note that the two Young-Laplace equations, the classical version (\ref{e:fb}) and the modified one (\ref{e:fb-v}), 
differ absolutely by 
\begin{equation}
\label{e:YL_diff}
\delta_{YL} = |p_* - p - (\gamma_* - \gamma) \mathcal H|.
\end{equation}
%$p^* - p$ and $\gamma^* - \gamma$, which are not small in general.
For a straight cylinder, with radius $r = 1$nm, we have at the interface $\partial_\tau \varphi = (0, 0, 0)^{t}$, as derived explicitly in \cite{BL-10}. 
Then 
${\displaystyle p_* - p = - \frac{1}{2} \frac{\sigma_c^2}{\varepsilon}}$ and %(with the 1st order approx. \varphi(R) \approx -RT/F \psi(x = 1) =
%RT/F * 2 \ln (1 - \lambda/8), 0 \leq \lambda < 8 from Berg/Ladipo 2009 for this radial symmetric situation)
$\gamma_* - \gamma = - 2 RT/F \, \sigma_c \ln(1- \lambda_p/8)$, where the parameter $\lambda_p = r^2/d_l^2$ depends on the so-called Debye length $d_l$ of the system.
Since $d_l > 0.6$nm \cite{PhysChemInterf}, we see that $0 \leq \lambda_p < 2.78$.
Thus (\ref{e:YL_diff}) reads 
\begin{equation*}
\delta_{YL} = \sigma_c \left| - \frac{1}{2} \frac{\sigma_c}{\varepsilon} + 2 \frac{RT}{F} %c_0
 \ln\left(1 - \frac{\lambda_p}{8} \right) \, \frac{1}{r} \right|.
\end{equation*}
Evidently, for small $\lambda_p$ we find that $\gamma^* - \gamma\approx -\frac{1}{2}\frac{\sigma_c^2}{\epsilon}$ is negligible, and then we expect similar solution with YL and modified YL equations. 
%\textcolor{blue}{(can we state the following comment by using the terms in the eq. above, and then eventually state that the difference of solutions can be large?)}
The relative difference for the Cauchy stress tensor is approximately $\delta_{YL}/{k_S}$ and is between $0.00356 $ and $0.00438 $. %$0.356 \%$ and $0.438 \%$. 
%$\lambda_p \to 8$ it would yield a singularity. %The paper says \lambda_p \approx 6.13, but that had been determined for a constant \eps. For the pressure, we find $p^* - p \approx 180 bar$, what is of the order of $\overline{p}$. 
As pointed out in \cite{LBKN-11}, we may expect some slight corrections to our numerical results by using a non-constant $\varepsilon_r$.


\noindent
\begin{thebibliography}{99}

\bibitem{adn-1+2}
{\sc S. Agmon, A. Douglis, L. Nirenberg}, 
\textit{Estimates near the boundary for solutions of
elliptic partial differential equations satisfying general boundary conditions}, 
Parts I and II, Comm. Pure Appl. Math. 12 %, no. 4
 (1959) 623-727; 17 %, no. 1
 (1964) 35-92.

\bibitem{bartolucci-1}
\textsc{D. Bartolucci, F. Leoni, L. Orsina, A.C. Ponce%Ponte
},
\textit{Semilinear equations with exponential nonlinearity and measure data},
Annales of Institute Henri Poincar\'e, AN 22 (2005) 779-815.

\bibitem{beatty-1}
\textsc{Millard F. Beatty},
\textit{A Lecture on Some Topics in Nonlinear Elasticity and Elasticity Stability},
Preprint Series \#99, Institute for Mathematics and its Applications, 
University of Minnesota, Minneapolis, 1984.
%\textcolor{blue}{(any of you knows a 'better' reference to refer to linear elasticity?)}

\bibitem{BL-10}
\textsc{P. Berg, K. Ladipo},
\textit{Exact solution of an electro-osmotic flow problem in a cylindrical channel of polymer electrolyte membranes},
Proc. Roy. Soc. A 465 (2009) 2663-2679.

\bibitem{brezis-1}
\textsc{H. Br\'ezis, W. Strauss},
\textit{Semi-linear second-order elliptic equations in $L^1$},
J . Math. Soc. Japan 25 (1973) 565-590. %, Vol. 25, No. 4, 1973

\bibitem{brezis-2}
\textsc{H. Br\'ezis, M. Marcus, A.C. Ponce},
\textit{Nonlinear elliptic equations with measures revisited},
%Chap. 4 in 
%\textit{Mathematical Aspects of Nonlinear Dispersive Equations (AM-163)},
%Edited by J.%ean
% Bourgain, C.%arlos E.
% Kenig \& S. Klainerman, Princeton, NJ, 2007
In: J. Bourgain, C. Kenig, S. Klainerman (Eds.), \textit{Mathematical Aspects of Nonlinear Dispersive Equations (AM-163)}, Princeton, 2007, pp. 55-110.


\bibitem{brezis-3}
\textsc{H. Br\'ezis},
\textit{Functional Analysis, Sobolev Spaces and Partial Differential Equations},
Springer, New York, 2010. %2011.

\bibitem{PhysChemInterf}
\textsc{H.-J. Butt, K. Graf, M. Kappl},
\textit{Physics and Chemistry of Interfaces},
Wiley, Weinheim, 2006.


\bibitem{ciarlet-1}
{\sc P.-G. Ciarlet},
{\it Mathematical Elasticity, Vol. 1, Three Dimensional Elasticity},
Elsevier, Amsterdam, 1988.

\colre{
\bibitem{ciarlet-2}
{\sc P.-G. Ciarlet},
{\it Lectures on Three-Dimensional Elasticity}, 
Tata Institute of Fundamental Research,
Bombay, 1983.
}



\bibitem{RecentTrendsCh3}
\textsc{K. S. Dhathathreyan, N. Rajalakshmi},
\textit{Polymer Electrolyte Membrane Fuel Cell},
In: S. %Suddhasatwa
Basu: Recent Trends in Fuel Cell Science and Technology, Springer, New York, 2007, pp. 40--115.

\bibitem{quarteroni-1}
\textsc{S. %imone
 Deparis, M. %arco
 Discacciati, G. %illes
 Fourestey, A. %lfio
 Quarteroni},
\textit{Fluid-structure algorithms based on Steklov-Poincar\'e operators},
Comput. Methods Appl. Mech. Engrg. 195 (2006) 5797-5812.

\bibitem{gilbarg+trudinger-1}
D. {\sc Gilbarg} and N.S. {\sc Trudinger},
{\it Elliptic Partial Differential Equations of Second Order},
{%\rm 
Springer, %-Verlag,
 Berlin, 1983}.
%\it <-> \rm

\bibitem{hunlich-1}
{\sc A. Glitzky, R. H\"unlich},
\textit{On energy estimates for electro-diffusion                                           
equations arising in semiconductor technology}, 
Research Notes in Mathematics 406 (2000) 158--174.


\bibitem{hiranao-1}
{\sc N. Hirano, W. Se Kim},
{\it Multiple existence of solutions for a semilinear elliptic problem with Neumann boundary condition},
%Journal of Mathematical Analysis and Applications,
%Volume 314, Issue 1, 1 February 2006, Pages 210-218 
J. Math. Anal. Appl. 314 (2006) 210-218.

\bibitem{KBN-12}
\textsc{S.-J. Kimmerle, P. Berg, A. Novruzi},
\textit{An electrohydrodynamic equilibrium shape problem for polymer electrolyte membranes in fuel cells},
In: System Modeling and Optimization %subtitle
 - 25th IFIP TC 7 Conference, Berlin, Germany, September 12-16, 2011, Revised Selected Papers, Approx. IX, 575 pp.,
IFIP AICT 391, Springer, Heidelberg, 2013.


\bibitem{stampacchia-1}
\textsc{D. Kinderlehrer, G. Stampacchia},
\textit{An Introduction to Variational Inequalities and Their Applications},
Academic Press, New York, 1980.


\bibitem{kreuer}
{\sc K.-D. Kreuer, S. Paddison, E. Spohr},
{\it Transport in proton conductors for fuel-cell applications: simulations,
elementary reactions, and phenomenology}, 
Chem. Rev. 104 (2004) 4637-4678.


\bibitem{LBKN-11}
{\sc K. Ladipo, P. Berg, S.-J. Kimmerle, A. Novruzi},
{\it Effects of radially-dependent parameters on proton transport in polymer electrolyte nanopores}, 
J. Chem. Phys. 134 (2011) 074103-1-12.



\bibitem{MGCR}
\textsc{E. Marchandise, P. Geuzaine, N. Chevaugeon, J.-F. Remacle},
\textit{A stabilized finite element method using a discontinuous level set approach for the computation of bubble dynamics},
J. Comp. Phys. %Journal of Computational Physics
 225 %(1),
(2006) 949--974.


\bibitem{mauritz}
{\sc K. Mauritz, R. Moore},
{\it State of understanding Nafion}, 
Chem. Rev. 104 (2004) 4535-4585.


\bibitem{SchmidtRohrChen}
\textsc{K. Schmidt-Rohr, Q. Chen},
\textit{Parallel cylindrical water nanochannels in Nafion fuel-cell membranes},
Nat. Mater. 7 (2008) 75--83.


\bibitem{SchmicklerSantos}
\textsc{W. Schmickler, E. Santos},
\textit{Interfacial Electrochemistry},
second ed., Springer, Heidelberg, 2010.


\bibitem{simon-1}
{\sc J. Simon},
\emph{ Differentiation with respect to the domain in boundary value problems},
Numer. Funct. Anal. Optim. 2 %, no 7-8
 (1980) 649-687.


\bibitem{simon-2}
{\sc J. Simon},
\emph{ Optimum design for Neumann condition and for related boundary value conditions},
In: J.-P. Zol\'esio (Ed.): Boundary Control and Boundary Variations, %ed. J.-P. Zol\'esio,
Proceedings of the IFIP WG 7.2 Conference, Nice, France, June 10-13, 1987,
Lecture Notes in Control and Information %Informations,
Sciences, vol. 100, Springer, Berlin, 1988.


\bibitem{Wilmanski-08}
\textsc{K. Wilmanski},
\textit{Continuum Thermodynamics, Part I: Foundations},
Series in Advances in Mathematics for Applied Sciences, vol.~77,
World Scientific Publishing, Singapore, 2008.

\end{thebibliography}
\end{document}